\documentclass[hidelinks,11pt]{article}

\newcommand{\bol}{\boldsymbol}

\newcommand{\nephi}{\boldsymbol{\phi}}

\newcommand{\nex}{\boldsymbol{x}}

\newcommand{\ner}{\boldsymbol{r}}

\newcommand{\de}{\,\mathrm{d}}                               
\newcommand{\e}{\operatorname{e}}                               
                     
\newcommand{\inc}{\mathrm{inc}}   
   
\newcommand{\gui}{\mathrm{gui}}   
\newcommand{\ontext}{\quad\mbox{on}\quad}   
   
\newcommand{\andtext}{\quad\mbox{and}\quad}

\newcommand{\p}{\partial}

\newcommand{\real}{\mathrm{Re}\,}    
                                   
\newcommand{\imag}{\mathrm{Im}\,}

\newcommand{\lf}{\left}
\newcommand{\rg}{\right}

\newcommand{\R}{\mathbb{R}}       
\newcommand{\C}{\mathbb{C}}

\usepackage{multirow}
\usepackage{amsmath}
\usepackage{amsfonts}
\usepackage{amsmath}
\usepackage{amssymb}  
\usepackage{graphicx}
\usepackage{caption}
\usepackage{mathrsfs}
\usepackage{upgreek}
\usepackage{amsthm}
\usepackage{subfig}
\usepackage{booktabs}
\usepackage{authblk}
\usepackage{MnSymbol}

\usepackage{color}

\newtheorem{theorem}{Theorem}[section]

\newtheorem{remark}[theorem]{Remark}

\title{Windowed Green Function method for the Helmholtz equation in
  presence of multiply layered media}
\author[1]{Oscar P. Bruno\thanks{obruno@caltech.edu}}
\author[2]{Carlos P\'erez-Arancibia}
\affil[1]{\small{Computing \& Mathematical Sciences, California
    Institute of Technology.}}  \affil[2]{\small{Department of
    Mathematics, Massachusetts Institute of Technology.}}
\date{\today}

\begin{document}
\maketitle

\begin{abstract}
This paper presents a new methodology for the
  solution of problems of two- and three-dimensional acoustic
  scattering (and, in particular, two-dimensional electromagnetic
  scattering) by obstacles and defects in presence an arbitrary number
  of penetrable layers. Relying on use of certain slow-rise windowing
  functions, the proposed Windowed Green Function approach (WGF)
  efficiently evaluates oscillatory integrals over unbounded domains,
  with high accuracy, without recourse to the highly expensive
  Sommerfeld integrals that have typically been used to account for
  the effect of underlying planar multi-layer structures.  The
  proposed methodology, whose theoretical basis was presented in the
  recent contribution (SIAM J. Appl. Math. 76(5), p. 1871, 2016), is
  fast, accurate, flexible, and easy to implement. Our numerical
  experiments demonstrate that the numerical errors resulting from the
  proposed approach decrease faster than any negative power of the
  window size. In a number of examples considered in this paper the
  proposed method is up to thousands of times faster, for a given
  accuracy, than corresponding methods based on use of Sommerfeld
  integrals.

\end{abstract}

\section{Introduction}


This paper presents a new methodology for the solution of problems of
acoustic scattering by obstacles and defects in presence an arbitrary
number of penetrable layers in two and three-dimensional space;
naturally, the two-dimensional Helmholtz solvers also apply, by
mathematical analogy, to corresponding two-dimensional electromagnetic
scattering problems. 
\maketitle	 \noindent This ``Windowed Green Function'' (WGF) method,
whose theoretical basis was presented in the recent
contribution~\cite{Bruno2015windowed}, is based on use of smooth
windowing functions and integral kernels that can be expressed
directly in terms of the free-space Green function, and, importantly,
it \emph{does not require use of expensive Sommerfeld integrals}.  The
proposed methodology is fast, accurate, flexible, and easy to
implement. Our experiments demonstrate that, as predicted by theory,
the numerical errors resulting from the proposed approach decrease
faster than any negative power of the window size. In a number of
examples considered in this paper the proposed method is up to
thousands of times faster, for a given accuracy, than corresponding
methods based on use of Sommerfeld integrals.

The classical layer Green functions and associated Sommerfeld
integrals automatically enforce the relevant transmission conditions
on the unbounded flat surfaces and thus reduce the scattering problems
to integral equations on the obstacles and/or defects
(cf.~\cite{Michalski:2016dz,perez2017windowed}). The Sommerfeld
integrals amount to singular Fourier
integrals~\cite{Chew1995waves,Sommerfeld1909} whose evaluation is
generally quite challenging. A wide range of approaches have been
proposed for evaluation of these
quantities~\cite{Aksun:2009fn,Cai:2002vt,Cai:2000bl,Cui:1998fw,Cui:1999tb,Lindell:1984es,Oneil2014efficient,Paulus:2000vr,PerezArancibia:2014fg}
but, as is known, all of these methods entail significant
computational
costs~\cite{Cai:2002vt,Lindell:1984es,Michalski:2016dz,Oneil2014efficient}.

The WGF approach proceeds as follows. The integral equation
formulations of the scattering problems under consideration, which are
at first posed on the complete set of material interfaces (including
all unbounded interfaces), are then {\em smoothly truncated} to produce
an approximating integral-equation system posed over bounded
integration domains that include the surface defects and relatively
small portions of the flat interfaces. The integral-operator
truncation is effected by means of a certain {\em slow-rise smooth
  window function} which, importantly, gives rise to solution errors
which decrease faster than any negative power of the window size. In
practice the proposed solution method is up to thousands of times
faster, for a given accuracy, than corresponding
methods~\cite{PerezArancibia:2014fg} based on use of Sommerfeld
integrals; the speedups in evaluations of near fields are even more
significant, in view of the large computing times required for
evaluation of Sommerfeld integrals near the planar interface.

This paper is organized as follows. Section~\ref{sec:ml_prelim}
presents a description of the multilayer scattering
problem. Section~\ref{sec:integral_equation_ML} then presents two
types of direct multi-layer integral equations for a physical, which
can be obtained by means of a generalized version of Green's third
identity (which is itself derived in
Appendix~\ref{eq:green_multi_layer}). The windowed integral equations
are derived in Section~\ref{sec:MWGFM} and the corresponding
expressions for the field evaluation are presented in
Section~\ref{sec:near_field_complicated_NL}. Section~\ref{sec:numerical_examples_ML},
finally, presents a variety of numerical examples which demonstrate
the super-algebraic convergence and the efficiency of the proposed
approach.

\section{Preliminaries\label{sec:ml_prelim}}
\begin{figure}[h!]
  \centering \subfloat[\label{fig:flat_domain}]{
\includegraphics[scale=0.8]{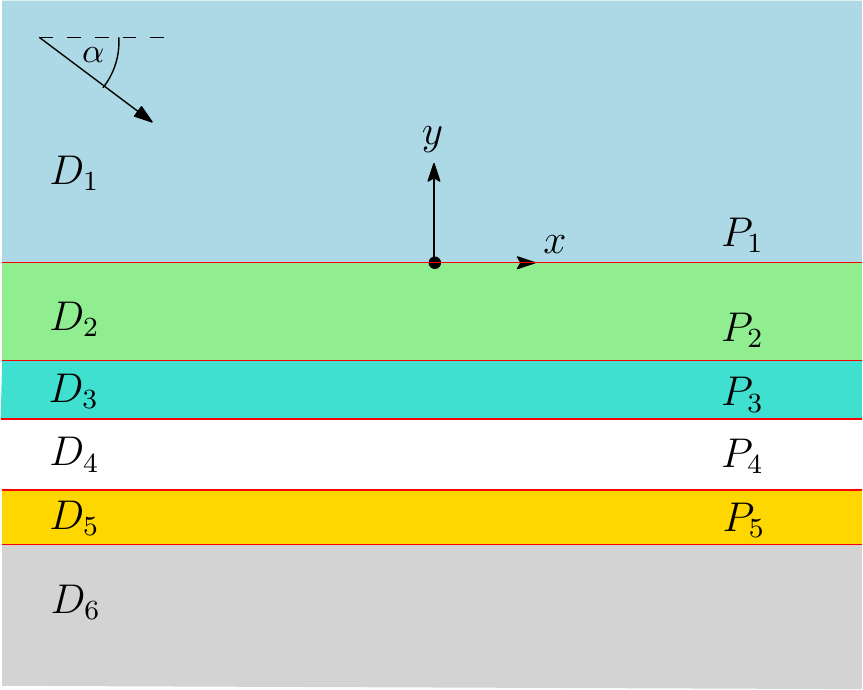}}\quad\subfloat[
	 \label{fig:domain}]{
           \includegraphics[scale=0.8]{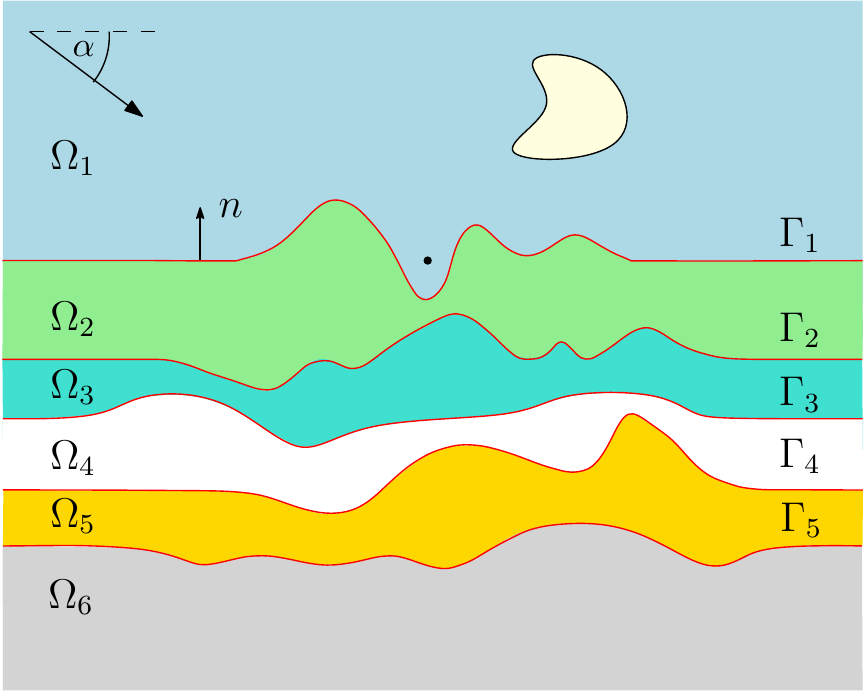}}
\caption{Geometry description of a two- or three-dimensional planar
  layered medium (a) and a locally perturbed planar layered medium (b)
  for the case $N=6$.}\label{fig:domains}
\end{figure}

We consider the problem of scattering of an acoustic incoming wave by
a two- or three-dimensional configuration such as the one depicted in
Figure~\ref{fig:domain}---in which an incoming wave is scattered by
localized (bounded) surface defects and/or scattering objects in
presence of a layered medium containing a number $N>1$ of layers. For
notational simplicity our descriptions are presented in the
two-dimensional case, but applications to three-dimensional
configurations are presented in
Section~\ref{sec:numerical_examples_ML}. The unperturbed
configuration, which is shown in Figure~\ref{fig:flat_domain} for
reference, consists of $N$ planar layers given by $D_j
=\R\times(-d_{j},-d_{j-1})$ for $j=1,\ldots, N$. The planar boundary
at the interface between the layers $D_j$ and $D_{j+1}$ is denoted by
$P_j=\R\times\{-d_j\}$ ($j=1,\ldots,N-1$). The corresponding perturbed
layers and their boundaries will be denoted by $\Omega_j$, $j=1,\ldots
N$ and $\Gamma_j$, $j=1,\ldots,N-1,$ respectively; naturally, it is
assumed that $\Gamma_j\cap\Gamma_i=\emptyset$.

Letting $\omega>0$ and $c_j>0$ denote the angular frequency and the
speed of sound in the layer $\Omega_j$, the wavenumber $k_j$ in that
layer is given by $k_j=\omega/c_j$.  Assuming e.g. an incident plane
wave of the form $u^{\inc}(\ner)=\e^{i k_1(x\cos\alpha +
  y\sin\alpha)}$ (where $\ner= (x,y)$ and where $-\pi < \alpha <0$
denotes the angle of incidence measured with respect to the
horizontal) and letting $u$ denote the acoustic pressure, the restrictions $u_j=u|_{\Omega_j}$
of the total field $u$ to the domains $\Omega_j$ ($j=1,\dots,N$)
satisfy the homogeneous Helmholtz equation
\begin{equation}
\Delta u_j+k_j^2u_j = 0\qquad\mbox{in}\qquad \Omega_j,\quad j=1,\dots,N,\label{eq:Helm_Eq}
\end{equation}
together with the transmission conditions 
\begin{equation}
 u_{j}= u_{j+1}\quad\mbox{and}\quad \frac{\p u_j}{\p n}=\nu_j\frac{\p u_{j+1}}{\p n}\quad\mbox{on}\quad\Gamma_j,\ \ j=1,\ldots,N-1,\label{eq:trans_conditions}
\end{equation}
where $\nu_j = \varrho_{j}/\varrho_{j+1}$ where $\varrho_j$ denotes the fluid density in $\Omega_j$. For definiteness,
here and throughout this paper the unit normal $n=n(\ner)$ for
$\ner\in\Gamma_j$ is assumed to point into $\Omega_j$.

As is well known, a closed form expression
exists~\cite{Chew1995waves,Stokes:1860vy} for the total field $u^p$
throughout space ($u^p = u^p_j$ in $D_j$, $j=1,\ldots,N$), that
results as a plane wave $u^\inc$ impinges on the planar layer medium
$D=\bigcup_{j=1}^N D_j$. In detail, letting $k_{1x}=k_1\cos\alpha$ and
$k_{jy}=\sqrt{k_j^2-k_{1x}^2}$, $j=1,\ldots, N$ (where the complex
square-root is defined in such a way that $\imag k_{jy}\geq 0$, which,
noting that $\imag k_j^2 \geq 0$, requires $\real k_{jy}\geq 0$ as
well), the planar-medium solution $u_j^p$ in $D_j$ is given by
\begin{equation}
u_j^p(x,y) = A_j\e^{ik_{1x}x}
\left\{\e^{-ik_{jy}y}+ \widetilde R_{j,j+1}\e^{ik_{jy}(y+2d_j)}\right\}\quad\mbox{in}\quad D_j,\ 1\leq j \leq N,
\label{eq:multi_layer_plane_wave}
\end{equation}
in terms of certain generalized reflection coefficients $\widetilde
R_{j,j+1}$ and amplitudes $A_j$. The amplitudes and the generalized
reflection coefficients can be obtained recursively by means of the
relations
\begin{eqnarray}
A_j&=&\left\{\begin{array}{cll}1&\mbox{if}&j=1,\medskip\\
\displaystyle\frac{T_{j-1,j}A_{j-1}\e^{i(k_{j-1,y}-k_{jy})d_{j-1}}}{1-R_{j,j-1}\widetilde R_{j,j+1}\e^{2ik_{jy}(d_j-d_{j-1})}}&\mbox{if}&j=2,\ldots,N,\end{array}\right.\label{eq:ampl}
\end{eqnarray}
and
\begin{eqnarray}
\widetilde R_{j,j+1}&=& \left\{\begin{array}{cll}0&\mbox{if}&j=N,\medskip\\
\displaystyle R_{j,j+1}+\frac{T_{j+1,j}\widetilde R_{j+1,j+2}T_{j,j+1}\e^{2ik_{j+1,y}(d_{j+1}-d_{j})}}{1-R_{j+1,j}\widetilde R_{j+1,j+2}\e^{2ik_{j+1,y}(d_{j+1}-d_{j})}}&\mbox{if}&j=N-1,\ldots,1,
\end{array}\right.\label{eq:gen_ref}
\end{eqnarray}
in terms of the reflection and transmission coefficients
\begin{eqnarray*}
R_{j,j+1} = \frac{k_{jy}-\nu_jk_{j+1,y}}{k_{jy}+\nu_jk_{j+1,y}}\quad\mbox{and}\quad T_{j,j+1} = \frac{2k_{jy}}{k_{jy}+\nu_jk_{j+1,y}},\label{eq:ref_trans_coef_TE}
\end{eqnarray*}
respectively.

\section{Integral equation
  formulations}\label{sec:integral_equation_ML}
This section presents an integral equation for the unknown interface
values of the total field and its normal derivative from below, at
each one of the interfaces $\Gamma_j$, $j=1,\ldots, N-1$. As in the
contribution~\cite{Bruno2015windowed} we utilize the single- and
double-layer potentials
\begin{equation}
\begin{split}
 S_j^{t}[\phi](\ner) = \int_{\Gamma_{j-1}} G_{k_j}(\ner,\ner')\phi(\ner')\de s_{\ner'},&\qquad
 D_j^{t}[\phi](\ner) = \int_{\Gamma_{j-1}} \frac{\p G_{k_j}}{\p n_{\ner'}}(\ner,\ner')\phi(\ner')\de s_{\ner'},\\
 S_j^{b}[\phi](\ner) = \int_{\Gamma_{j}} G_{k_j}(\ner,\ner')\phi(\ner')\de s_{\ner'},&\qquad
 D_j^{b}[\phi](\ner) = \int_{\Gamma_{j}} \frac{\p G_{k_j}}{\p n_{\ner'}}(\ner,\ner')\phi(\ner')\de s_{\ner'},
\end{split}\label{eq:layer_potentials_NL}
\end{equation}
which are defined for $\ner\in\R^2$ and are expressed in terms of improper integrals whose convergence
is conditioned upon the oscillatory behavior of the integrand.  Here we
have called
$G_{k_j}(\ner,\ner')=\frac{i}{4}H_0^{(1)}(k_j|\ner-\ner'|)$ the
free-space Green function for the Helmholtz equation with wavenumber
$k_j$. Additionally, we define the integral operators
\begin{equation}
\begin{split}
 K_j^{t}[\phi](\ner) = \int_{\Gamma_{j-1}} \frac{\p G_{k_j}}{\p n_{\ner}}(\ner,\ner')\phi(\ner')\de s_{\ner'}&\qquad
 N_j^{t}[\phi](\ner) = \int_{\Gamma_{j-1}} \frac{\p^2 G_{k_j}}{\p n_{\ner}\p n_{\ner'}}(\ner,\ner')\phi(\ner')\de s_{\ner'},\\
 K_j^{b}[\phi](\ner) = \int_{\Gamma_{j}} \frac{\p G_{k_j}}{\p n_{\ner}}(\ner,\ner')\phi(\ner')\de s_{\ner'}&\qquad
 N_j^{b}[\phi](\ner) = \int_{\Gamma_{j}} \frac{\p^2 G_{k_j}}{\p n_{\ner}\p n_{\ner'}}(\ner,\ner')\phi(\ner')\de s_{\ner'},
\end{split}\label{eq:diff_layer_potentials_NL}
\end{equation}
where the evaluation point $\ner$ belongs to either $\Gamma_j$ or $\Gamma_{j-1}$.

In order to formulate an integral equation for the unknown interface
values we define the unknown density functions $\varphi_j:\Gamma_j\to\C$ and
$\psi_j:\Gamma_j\to\C$ ($j=1,\ldots,N-1$) by
\begin{equation}
\varphi_j=u_{j+1}\quad \mbox{and}\quad \psi_j=\frac{\p u_{j+1}}{\p n}\quad\mbox{on}\quad \Gamma_j.\label{eq:unknowns}
\end{equation}
Additionally we define the vector density functions
\begin{equation}
\boldsymbol\phi_{j}=\left[ \varphi_j, \psi_j\right]^T,\quad
\boldsymbol\phi^{\mathrm{inc}}=\lf[ u^{\mathrm{inc}}\big|_{\Gamma_1}, \frac{\p
u^\inc}{\p n}\Big|_{\Gamma_1}\rg]^T\quad \mbox{and}\quad\bol\phi^{\parallel}=\lf[
u_N^{\parallel}\big|_{\Gamma_{N-1}},\frac{\p u^{\parallel}_N}{\p
n}\Big|_{\Gamma_{N-1}}\rg]^T\label{eq:vector_density}
\end{equation}
where $u_N^\parallel$ is defined in~\eqref{eq:parallel_ML}, and the matrix operators
\begin{equation}
\begin{split}
&\operatorname E_{j}=\left[\begin{array}{cc}
1  & 0 \\
0&\frac{1+\nu_j}{2}
\end{array}\right],\;\;\;\; \operatorname T_{j}=\left[\begin{array}{cc}
D_{j+1}^{t} - D_j^{b} &  - S_{j+1}^{t}+\nu_j S_j^{b} \smallskip\\
N_{j+1}^{t} - N_j^{b} &- K_{j+1}^{t}+\nu_j K_j^{b}  
\end{array}\right],\ \operatorname L_{j}=\left[\begin{array}{cc}
D_{j}^{t}  & -S_{j}^{t}\smallskip \\
N_{j}^{t}& -K_{j}^{t}
\end{array}\right]\\
&\mbox{and }\operatorname R_{j}=\left[\begin{array}{cc}
-D_{j+1}^{b}  & \nu_{j+1}S_{j+1}^{b} \smallskip\\
-N_{j+1}^{b}& \nu_{j+1}K_{j+1}^{b}
\end{array}\right]\mbox{ (all operators evaluated at observation points $\ner$ on $\Gamma_j$).}
\end{split}
\label{eq:notation_system}\hspace{-.4cm}
\end{equation}

A general multi-layer integral formulation of the
problem~\eqref{eq:Helm_Eq}--\eqref{eq:trans_conditions} can now be
obtained in terms of these densities and operators. Indeed, as is
shown in Appendix A, the fields within the layers admit the integral
representations
\begin{equation}
\begin{array}{lclll}
  u_1(\ner) &=&  D_{1}^{b}[\varphi_1](\ner)-\nu_1 S_{1}^{b}[\psi_1] (\ner)+u^{\inc}(\ner),&\smallskip\\
  u_j(\ner) &=&  D_{j}^{b}[\varphi_j](\ner)-\nu_j S_{j}^{b}[\psi_j] (\ner)\smallskip\\
  &&- D_{j}^{t}[\varphi_{j-1}](\ner)+ S_{j}^{t}[\psi_{j-1}] (\ner),& j=2,\ldots,N-1,\smallskip\\
  u_N(\ner) &=&  - D_{N}^{t}[\varphi_{N-1}](\ner)+ S_{N}^{t}[\psi_{N-1}](\ner)+u^\parallel_N(\ner),&\end{array}\label{eq:rep_ML}\end{equation}
in terms of the interface values~\eqref{eq:unknowns}. Therefore,
evaluating $u_1 + u_2$ and $\p u_1/\p n + \p
u_2/\p n$ on $\Gamma_1$ from the boundary values on $\Gamma_1$ of the
expressions in~\eqref{eq:rep_ML} and their
normal derivatives, and using the notations~\eqref{eq:vector_density}
and~\eqref{eq:notation_system}, we obtain the $j=1$ interface equation
\begin{subequations}
\begin{equation}\label{int_eq_1}
  \operatorname E_1 \boldsymbol \phi_1 + \operatorname T_{1}\left[\boldsymbol\phi_1\right] + \operatorname R_{1}\left[\boldsymbol\phi_2\right] = \boldsymbol\phi^\inc\ontext \Gamma_1.
\end{equation}
A similar procedure
yields the integral equations
\begin{equation}
\operatorname E_j \boldsymbol\phi_j +\operatorname L_{j}\left[\boldsymbol\phi_{j-1}\right]+ \operatorname T_{j}\left[\boldsymbol\phi_j\right] + \operatorname R_{j}\left[\boldsymbol\phi_{j+1}\right] =\bol 0\ontext\Gamma_j,\quad j=2,\ldots,N-2
\end{equation} 
and
\begin{equation}
  \operatorname E_{N-1} \boldsymbol\phi_{N-1} + \operatorname L_{N-1}\left[\boldsymbol\phi_{N-2}\right] + \operatorname T_{N-1}\left[\boldsymbol\phi_{N-1}\right] = \nephi^\parallel\ontext\Gamma_{N-1}. \end{equation}\label{eq:separate_layers}\end{subequations}
(Note that, of course, the calculations leading to  equations~\eqref{eq:separate_layers} rely on the
well-known jump relations for the single- and double-layer potentials
and their normal derivatives~\cite{COLTON:1983}.)  

\begin{remark}\label{multi_vectors}
  In what follows equations~\eqref{eq:separate_layers} are expressed
  in terms of a single column vector function $\bol\phi$ (defined on
  the Cartesian product $\Gamma=\prod_{j=1}^{N-1}\Gamma_j$ of the
  curves $\Gamma_j$) whose $j$-entry equals the density pair
  $\boldsymbol\phi_{j}=[\varphi_j,\psi_j]^T:\Gamma_j\to\C^2$ for
  $j=1,\dots, N-1$. We may thus write
  $$\boldsymbol\phi = [\phi_{1},\psi_1,\varphi_2,\psi_2,\cdots
  ,\varphi_{N-1},\psi_{N-1}]^T:\Gamma\to\C^{2(N-1)}.$$ Similarly we define
$$\boldsymbol\phi^\inc =
  \lf[u^\inc|_{\Gamma_1},\frac{\p u^\inc}{\p n}\Big|_{\Gamma_1},0,0,\cdots
 ,0,0 ,u^\parallel_{N}|_{\Gamma_{N-1}},\frac{\p u^\parallel_N}{\p n}\Big|_{\Gamma_{N-1}}\rg]^T:\Gamma\to\C^{2(N-1)}.$$ 

 With a slight notational abuse we will write $\bol\phi
 =[\bol\phi_1,\bol\phi_2,\ldots,\bol\phi_{N-1}]^T =
 [\phi_{1},\psi_1,\varphi_2,\psi_2,\cdots
 ,\varphi_{N-1},\psi_{N-1}]^T$. More generally, given arbitrary
 vectors $\bol\mu_j=[\alpha_j,\beta_j]^T:\Gamma_j\to\C^2$ for
 $j=1,\dots, N-1$ we will use the ``block-vector'' notation $\bol\mu
 =[\bol\mu_1,\bol\mu_2,\ldots,\bol\mu_{N-1}]^T =
 [\alpha_{1},\beta_1,\alpha_2,\beta_2,\cdots
 ,\alpha_{N-1},\beta_{N-1}]^T:\Gamma\to\C^{2(N-1)}$.

\end{remark}

Using the operators 
\begin{equation}
\mathcal E = \left[\begin{array}{cccc}
\operatorname E_1&&&\\
&\operatorname E_2&&\\
&&\ddots&\\
&&&\operatorname E_{N-1}
\end{array}\right]
\andtext \mathcal T_\Gamma = \left[\begin{array}{ccccc}
\operatorname T_1&\operatorname R_1&&&\\
\operatorname L_2&\operatorname T_2&\operatorname R_2&&\\
&\operatorname L_3&\ddots&\ddots\\
&&\ddots&\ddots&\operatorname R_{N-2}\\
&&&\operatorname L_{N-1}&\operatorname T_{N-1}
\end{array}\right]\label{eq:Tbig}\end{equation}
together with the notations~introduced in Remark~\ref{multi_vectors},
equations~\eqref{eq:separate_layers} can be expressed in the form 
\begin{equation}
\mathcal E \boldsymbol\phi+\mathcal T_\Gamma\left[\boldsymbol\phi\right] =\boldsymbol{\phi}^{\mathrm{inc}}\ontext \Gamma.\label{eq:transmission_system}
\end{equation}

\section{Windowed integral equations}\label{sec:MWGFM}
Following~\cite{Bruno2015windowed}, in this section we introduce
rapidly-convergent windowed versions of the integral
formulation~\eqref{eq:transmission_system}. In order to do so we utilize
the $(N-1)\times (N-1)$ block-diagonal matrix-valued window function
$\mathcal W_A:\prod_{j=1}^{N-1}\Gamma_{j}\mapsto \R^{2(N-1)\times 2(N-1)}$ given by
\begin{equation}
\mathcal W_A(\ner_1,\ner_2,\ldots,\ner_{N-1})  = \left[\begin{array}{cccc}
w_A(x_1)\operatorname I&&&\\
&w_A(x_2)\operatorname I&&\\
&&\ddots&\\
&&&w_A(x_{N-1})\operatorname I
\end{array}\right],\quad \ner_j=(x_j,y_j)\in\Gamma_j,\label{eq:window_matrix_NL}
\end{equation}
 in terms of the two-by-two identity matrix $\operatorname I$ and the
smooth window function
\begin{equation}\label{window_fnct}
  w_A(x) = \eta(x/A ;c,1),
\end{equation}
where $0<c<1$ and where
\begin{equation}
\eta(t;t_0,t_1) =\left\{
\begin{array}{cll}
1,&|t|\leq t_0,\medskip\\
\!\!\exp\left(\displaystyle\frac{2\e^{-1/u}}{u-1}\right),&\displaystyle t_0<|t|<t_1, u=\frac{|t|-t_0}{t_1-t_0},\medskip\\
0,&|t|>t_1.
\end{array}\right.\label{eq:window_function}
\end{equation}

Clearly $\eta$ and $w_A$ are infinitely differentiable
compactly-supported functions of $x$ and $t$, respectively. The
support of the window function $w_A = w_A(x)$ as a function of $\ner =
(x,y)\in \R^2$ equals the set $[-A,A]\times \R$. Note that the parameter $c$,
which controls the steepness of the rise of the window function $w_A$,
is not displayed as part of the notation $w_A$.

(While different values $A_j$ of the window-size $A$ ,
$j=1,\dots,N-1$, could in principle be used for the various layer
interfaces and corresponding block entries
in~\eqref{eq:window_matrix_NL}---possibly utilizing smaller
(resp. larger) values $A_j$ in higher (resp. lower) frequency layers,
and therefore reducing the overall number of unknowns required for the
WGF method to produce a given accuracy.  For simplicity, however,
throughout this paper a single window-size value $A$ is used for all
the interfaces.)

In order to produce a windowed version of
equation~\eqref{eq:transmission_system} we proceed in two stages. At first the integrand
is multiplied by the window matrix $\mathcal W_A$ and the equation is
restricted to the windowed region $\Gamma_A = \{(
\ner_1,\dots,\ner_{N-1})\in\Gamma:
\ner_j=(x_j,y_j)\mbox{ and }w_A(x_j)\ne 0 \mbox{ for all }j\}\subset
\R^{2(N-1)}$---so that, moving the remainder of the windowed integral
operator to the right-hand side and letting $\mathcal  I$ denote the
identity matrix of dimension $2(N-1)\times 2(N-1)$, the exact relation
\begin{equation}
\mathcal E \bol\phi+\mathcal T_\Gamma\left[\mathcal W_A\bol\phi\right] =\bol{\phi}^{\inc}-\mathcal T_\Gamma\left[( \mathcal I-\mathcal W_A)\bol\phi\right]\ontext\Gamma_A\label{eq:transmission_system_window_prev}
\end{equation} 
results. Note that defining 
$\Gamma_{j,A} = \Gamma_j\cap \{w_A\neq
0\}=\Gamma_j\cap\left\{ [-A,A]\times \R\right\}$
we have
\begin{equation}\label{gamma_A}
\Gamma_A = \prod_{j=1}^{N-1}\Gamma_{j,A}.
\end{equation} 

A successful implementation of the WGF idea requires use of an
accurate substitute for the quantity $\mathcal T_\Gamma\left[(
  \mathcal I-\mathcal W_A)\bol\phi\right]$ throughout $\Gamma_A$,
which does not depend on knowledge of the unknown density $\bol\phi$
(cf.~\cite{Bruno2015windowed}). In order to obtain such an
approximation we introduce an operator $\mathcal T_P$ which is
defined just like $\mathcal T_\Gamma$ in~\eqref{eq:Tbig} but in terms
of potentials~\eqref{eq:layer_potentials_NL} and
operators~\eqref{eq:diff_layer_potentials_NL} given by integrals on
the flat interfaces $P_j$ depicted in
Figure~\ref{fig:flat_domain}. Since $( \mathcal I-\mathcal W_A)$
vanishes wherever $\Gamma$ differs from $P=\prod_{j=1}^{N-1}P_j$, we clearly have
$\mathcal T_\Gamma\left[( \mathcal I-\mathcal
  W_A)\bol\phi\right]=\mathcal T_P\left[( \mathcal I-\mathcal
  W_A)\bol\phi\right]$.  Additionally, we consider the
aforementioned scalar densities $\varphi^p_j=u^p_{j+1}|_{\Gamma_j}$
and $\psi^p_j=\p u^p_{j+1}/\p n|_{\Gamma_j}$ on $P_j$
($j=1,\ldots,N-1$) that are associated with the planarly layered
medium $P$.  As shown in~\cite{Bruno2015windowed}, letting
$[\boldsymbol\phi^p]_j=[\varphi^p_j,\psi^p_j]^T$ ($j=1,\ldots,N-1$),
substitution of $\bol\phi$ by $\bol\phi^p$ on the right-hand-side
of~\eqref{eq:transmission_system_window_prev} results in errors that
decay super-algebraically fast as $A\to \infty$
within the subset
\begin{equation}\label{w_eq_1}
  \widetilde\Gamma_A = \Gamma_A\cap\prod_{j=1}^{N-1}\{(x_j,y_j)\in\Gamma_j:w_A(x_j)=1\}
\end{equation}
of $\Gamma_A$ wherein the window function $w_A$ equals one. Indeed,
even though $\bol\phi$ may differ significantly from $\bol\phi^p$, the
corresponding integrated terms result in super-algebraically small
errors, as it may be checked via stationary phase analysis
(see~\cite{Bruno2015windowed}, for details). We thus obtain the
super-algebraically-accurate windowed integral equation system
\begin{equation}
\mathcal E \bol\phi^w+\mathcal T_\Gamma\left[\mathcal W_A\bol\phi^w\right] =\bol{\phi}^{\inc}-\mathcal T_P[( \mathcal I-\mathcal W_A)\bol\phi^p]\ontext \Gamma_A
\label{eq:transmission_system_window_prev}
\end{equation} 
which
we re-express in the form
\begin{equation}
\mathcal E \bol\phi^w+\mathcal T_\Gamma\left[\mathcal W_A\bol\phi^w\right] =\bol{\phi}^{\inc} +  \mathcal T_P\left[\mathcal W_A\bol\phi^p\right]-\mathcal T_P\left[\bol\phi^p\right]\ontext \Gamma_A.\label{eq:transmission_system_window}
\end{equation} 
As shown in what follows, the right-hand term $\mathcal
T_P\left[\bol\phi^p\right]$
in~\eqref{eq:transmission_system_window} can be expressed in closed
form, and thus, using numerical integration over the bounded domain
$\Gamma_A$ to produce the term $\mathcal T_P\left[\mathcal
  W_A\bol\phi^p\right]$, the complete right-hand side can be
efficiently evaluated for any given $\nex\in\Gamma_A$.

A closed-form expression for $\bol\mu = \mathcal
T_P\left[\bol\phi^p\right]$ (cf. Remark~\ref{multi_vectors}) can be
obtained via an application of Green's formula: using~\eqref{Green-f}
with $C=P_j$ ($j=1,\dots,N-1$),
equations~\eqref{eq:green_representation_NL} yield the desired
relations:
\begin{equation}
    \bol\mu_j =\delta_{1,j}\nephi^\inc+\delta_{N-1,j}\nephi^\parallel-\left\{\begin{array}{clc} \operatorname E_j \bol\phi^p_j&\mbox{on}\quad\Gamma_j\cap P_j,\medskip\\
        \left[\begin{array}{cc}u^p\\  \nabla u^p\cdot n\end{array}\right]&\mbox{on}\quad\Gamma_j\cap(D_j\cup D_{j+1}),
 \end{array}\right.
\label{eq:correction_ML}
\end{equation}
for $j=1,\ldots,N-1,$ where $\delta_{i,j}$ denotes the Kronecker delta symbol.

As demonstrated in Section~\ref{sec:numerical_examples_ML} through a
variety of numerical examples, the vector density function
$\bol\phi^w$, which is the solution of the windowed integral
equation~\eqref{eq:transmission_system_window}, converges
super-algebraically fast to the exact solution $\bol\phi$
of~\eqref{eq:transmission_system} within $\Gamma^1_A$ as
the window size $A>0$ increases.  This observation can be justified
via arguments analogous to those presented in~\cite{Bruno2015windowed}.

\begin{remark}
  The difference $N^{t}_{j+1}-N^{b}_j$ of hypersingular operators that
  appears in the definition of the diagonal blocks $\operatorname T_j
  $ of~$\mathcal T_\Gamma$ is in fact a weakly singular integral
  operators (cf.~\cite[Sec. 3.8]{COLTON:1983}).
\end{remark}

\section{Near-field evaluation}\label{sec:near_field_complicated_NL}

This section presents a super-algebraically accurate WGF approximation
$u^w$ of the solution $u$
of~\eqref{eq:Helm_Eq}--\eqref{eq:trans_conditions} near the localized
defects. In order to obtain this approximation we consider the
``defect'' field
\begin{equation}
u^d_j=u_j-\tilde u^p_j\quad\mbox{in}\quad \Omega_j\quad (j=1,\ldots,N-1),\label{eq:def_field}
\end{equation} 
given by
the difference between the total field $u_j$ and the planar-structure
total field
\begin{equation}
  \tilde u_j^p(x,y) = A_j\e^{ik_{1x}x}
  \left\{\e^{-ik_{jy}y}+ \widetilde R_{j,j+1}\e^{ik_{jy}(y+2d_j)}\right\}\quad\mbox{in}\quad\Omega_j,\ 1\leq j \leq N.
\label{eq:multi_layer_plane_wave_2}
\end{equation}
Note that $\tilde u_j^p$ is given in $\Omega_j$ by the expressions on
the right-hand side of equation~\eqref{eq:multi_layer_plane_wave}.

Subtracting the integral representation
\begin{equation}
\begin{array}{lcll}
 \tilde u^p_1(\ner) &=&   D_{1}^{b}\lf[\tilde\varphi^p_1+f_1\rg](\ner)-\nu_1  S_{1}^{b}\lf[\tilde\psi^p_1+g_1\rg] (\ner)+u^{\inc}(\ner),&\smallskip\\
\tilde u^p_j(\ner) &=&   D_{j}^{b}\lf[\tilde\varphi^p_j+f_j\rg](\ner)-\nu_j  S_{j}^{b}\lf[\tilde\psi^p_j+g_j\rg] (\ner)\smallskip\\ 
&&-  D_{j}^{t}\lf[\tilde\varphi^p_{j-1}\rg](\ner)+  S_{j}^{t}\lf[\tilde\psi^p_{j-1}\rg] (\ner),&   j=2,\ldots,N-1,\smallskip\\
\tilde u^p_N(\ner)&=& -  D_{N}^{t}\lf[\tilde\varphi^p_{N-1}\rg](\ner)+  S_{N}^{t}\lf[\tilde\psi^p_{N-1}\rg](\ner)+u^\parallel_N(\ner),&
\end{array} \label{eq:rep_ML_alternative}\end{equation}
---which follows as equation~\eqref{eq:green_representation_NL} is applied to $\tilde
u^p_j$---from the integral representation~\eqref{eq:rep_ML}
we obtain the exact integral relations
    \begin{equation}\begin{array}{lcll}
        u^d_1(\ner) &=&  D_{1}^{b}[\varphi_1-\tilde\varphi^p_1-f_1](\ner)-\nu_1 S_{1}^{b}[\psi_1-\tilde\psi^p_1-g_1] (\ner),&\smallskip\\
        u^d_j(\ner)  &=&  D_{j}^{b}[\varphi_j-\tilde\varphi^p_j-f_j](\ner)-\nu_j S_{j}^{b}[\psi_j-\tilde\psi^p_j-g_j] (\ner)&\smallskip\\
        &&- D_{j}^{t}[\varphi_{j-1}-\tilde\varphi^p_{j-1}](\ner)+ S_{j}^{t}[\psi_{j-1}-\tilde\psi^p_{j-1}] (\ner),   
        & j=2,\ldots,N-1,\smallskip\\
        u^d_N(\ner) &=&  - D_{N}^{t}[\varphi_{N-1}-\tilde\varphi^p_{N-1}](\ner)+ S_{N}^{t}[\psi_{N-1}-\tilde\psi^p_{N-1}](\ner)&\end{array}\label{eq:surb_prelim}\end{equation}
    for the defect fields. These relations can be used to evaluate the defect fields $u^d_j$ in terms of the solution~$\bol\nephi_j=[\varphi_j,\psi_j]^T$ of the integral equation~\eqref{eq:transmission_system} together with the planar-structure total fields 
\begin{equation}\label{phi-psi_tilde}
  \tilde{\bol\phi}^p_j
=\left[\tilde\varphi^p_j,\tilde\psi^p_j\right]^T,\quad\mbox{where}\quad\tilde \varphi_j^p=\tilde u_{j+1}^p|_{\Gamma_{j}} \quad \mbox{and}\quad\tilde\psi_j^p = \frac{\p\tilde u^p_{j+1}}{\p n}\Big|_{\Gamma_j},
\end{equation}
and the jumps
\begin{equation}
\bol\uppsi_j =\left[f_j , g_j\right]^T,\quad\mbox{where}\quad f_j = \tilde u^p_{j}-\tilde u^p_{j+1}\andtext g_j =\frac{1}{\nu_j}\frac{\p \tilde u^p_j}{\p n}-\frac{\p \tilde u^p_{j+1}}{\p n}\quad\mbox{on}\quad \Gamma_j.\label{eq:defintion_jumps}
\end{equation}  
Note that, importantly, for each $j$, $1\leq j\leq N$, the functions
$f_j$ and $g_j$ vanish outside the $j$-th portion
$\Gamma_j\setminus\Pi_j$ of the boundary of the localized defects.

Relying on the WGF solutions $\boldsymbol\phi^{w}$ of
equation~\eqref{eq:transmission_system_window} and applying a
windowing procedure similar to the one used in the previous section, a
highly-accurate approximation to the defect
near-fields~\eqref{eq:surb_prelim} results. In detail, substitution of
$\varphi_j$ by $w_A\varphi_j^{w}+(1-w_A)\tilde\varphi^p_j$ and
$\psi_j$ by $w_A\psi_j^{w}+(1-w_A)\tilde\psi^p_j$
in~\eqref{eq:surb_prelim} yields the approximate expressions
\begin{equation}
\begin{array}{lcl}
u^{d,w}_1(\ner) &=& D_1^{b}\left[w_A(\varphi^{w}_1-\tilde\varphi^p_{1})-f_1\right]-\nu_1 S_1^{b}\left[w_A(\psi^{w}_1-\tilde\psi^p_{1})-g_1\right], \smallskip\\
u^{d,w}_j(\ner)&=&   D_{j}^{b}\lf[w_A(\varphi^{w}_j-\tilde\varphi^p_{j})-f_j\rg]-\nu_j S_{j}^{b}\lf[w_A(\psi^{w}_j-\tilde\psi^p_{j})-g_j\rg]\smallskip \\
&&- D_{j}^{t}\lf[w_A(\varphi^{w}_{j-1}-\tilde\varphi^p_{j-1})\rg]+ S_{j}^{t}\lf[w_A(\psi^{w}_{j-1}-\tilde\psi^p_{j-1})\rg],\quad j=2,\ldots,N-1,\smallskip\\
u_N^{d,w}(\ner)&=&- D_N^{t}\left[w_A(\varphi^{w}_{N-1}-\tilde\varphi^p_{N-1})\right]+ S_N^{t}\left[w_A(\psi^{w}_{N-1}-\tilde\psi^p_{N-1})\right],
\end{array}\label{eq:approx_defect_field}
\end{equation}
for the defect field $u^d_j$. The desired approximation $u^w_j$ for
the total field $u_j$ then follows from~\eqref{eq:def_field}:
 \begin{equation} 
u^w_j = \tilde u^p_j +
    u^{d,w}_j\quad\mbox{in}\quad \Omega_j\quad (j=1,\ldots,N).\label{eq:simpler_near_field_evaluation_NL}
\end{equation}

Formulae~\eqref{eq:simpler_near_field_evaluation_NL} provide
super-algebraically accurate approximations of the total near-fields
within the region
\begin{equation}\label{omega_eq_1}
  \widetilde\Omega_A = \bigcup_{j=1}^{N}\Omega_j\cap\{\ner\in\R^2:w_A(x)=1\}
\end{equation}
containing the localized defects---with uniformly small errors, as
$A\to\infty$, within every bounded subset of $\widetilde\Omega_A$. A
theoretical discussion in these regards (for the two-layer case) can
be found in~\cite{Bruno2015windowed} (see e.g. Remark~4.1 in that
reference).

\section{Far-field evaluation}\label{sec:far_field_complicated_NL}

\begin{figure}[h!]
\centering	
	\includegraphics[scale=0.8]{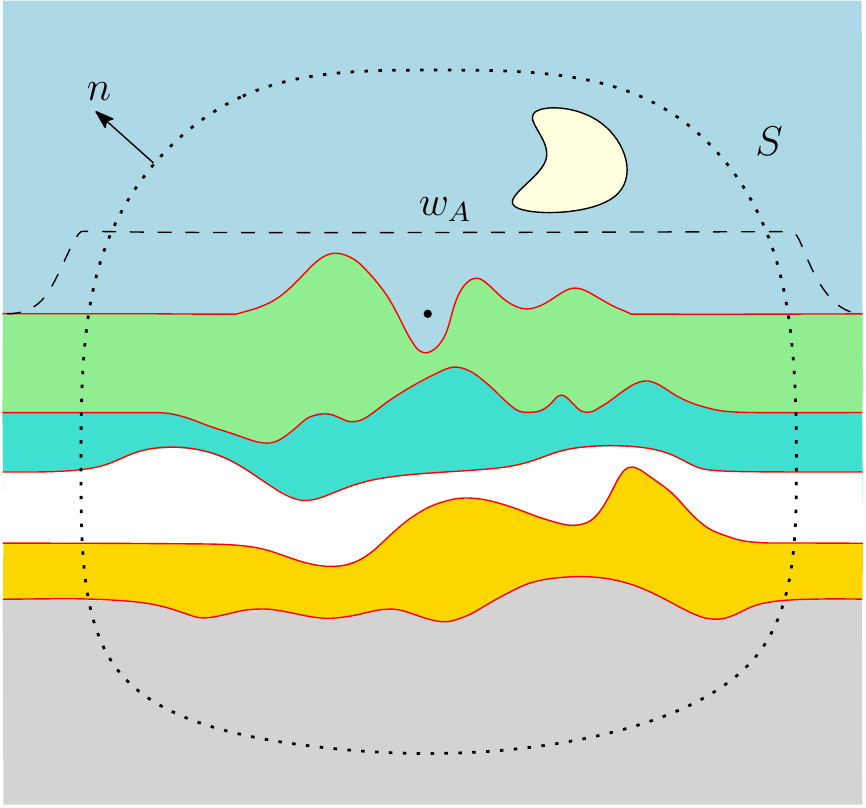}
\caption{Curve $S$  in~\eqref{eq:ext_field}.}\label{fig:ext_domain}
\end{figure}

As indicated in the previous section,
formulae~\eqref{eq:approx_defect_field}--\eqref{eq:simpler_near_field_evaluation_NL}
only provide uniformly accurate approximations within bounded subsets
of $\widetilde\Omega_A$. But, once accurate defect fields $u^{d,w}_j$
($j=1,\ldots,N$) have been obtained within $\widetilde\Omega_A$,
correspondingly accurate far-field values for the solution $u$ can be
obtained by applying certain Green-type formulae on a bounding curve
$S$, such as the one depicted in Figure~\ref{fig:ext_domain}, which
encloses all of the local defects, and which is contained within
$\widetilde\Omega_A$. In detail, defining the defect field
$u^d=u^d(\ner)$ to equal $u^d_j(\ner)$ for $\ner\in \Omega_j$
($j=1,\ldots,N$), use of a Green identity based on the $N$-layer Green
function~$H$ over the region exterior to $S$ leads to the integral
representation~\cite[Lemma 4.2.6]{perez2017windowed}
\begin{equation}
  u^{d}(\ner) = \int_{S}\lf\{\frac{\p H}{\p n_{\ner'}}(\ner,\ner')u^{d}(\ner')-H(\ner,\ner')\frac{\p u^{d}}{\p n}(\ner')\rg\}\de s_{\ner'},\label{eq:ext_field}
\end{equation}
which is valid for $\ner$ everywhere outside $S$.  Note that the
necessary values of $u^d_j$ and their normal derivatives on $S$ can be
computed by means
of~\eqref{eq:simpler_near_field_evaluation_NL}---since, by
construction, $S$ lies inside the region
where~\eqref{eq:simpler_near_field_evaluation_NL} provides an accurate
approximation of the field $u^d_j$.

The far-field approximation $u^f$ of the defect field $u^d$ as
$\ner\to\infty$ in any direction is then obtained by replacing the
layer Green function $H$ and its normal derivative $\p H/\p n_{\ner'}$
in~\eqref{eq:ext_field} by the respective first-order
$|\ner|\to\infty$ asymptotic expansions $H^f$ and $\p H^f/\p
n_{\ner'}$---which can be obtained for the $N$-layer case (as
illustrated
in~\cite{Bleistein1975Asymptotic,Chew1995waves,perez2017windowed,Bruno2015windowed,Brekhovskikh2013Acoustics}
for $N=2$ and below in this section for $N=3$) by means of the method
of steepest descents.  (The fact that the far field of the function
$\p H/\p n_{\ner'}$ coincides with $\p H^f/\p n_{\ner'}$ can be
verified by direct inspection of these two quantities.) The far field
$u^f$ is thus given by
\begin{equation}
  u^f(\ner) = \int_{S}\lf\{\frac{\p H^f}{\p n_{\ner'}}(\ner,\ner')u^{d}(\ner')-H^f(\ner,\ner')\frac{\p u^{d}}{\p n}(\ner')\rg\}\de s_{\ner'}.\label{eq:ext_field_ff}
\end{equation}
It is important to note that, unlike the layer Green function $H$
itself, the corresponding far-field $H^f$ and its normal derivative
can be evaluated inexpensively by means of explicit expressions. 

As an example we sketch here the calculation of the far-field $H^f$
for a slab---that is, a three-layer medium with wavenumbers $k_j$,
$j=1,2,3$ where $k_1=k_3$---in two-dimensional space. We assume the
case $k_2>k_1$ for which the slab can sustain guided modes that
propagate along the $x$-axis. In order to evaluate $H^f$ we first note
that, for a source point $\ner'=(x',y')\in D_j$ and a target point
$\ner=(r\cos\theta,r\sin\theta)\in D_1$ ($\theta\in[0,\pi]$), the
layer Green function $H$ is given by the contour
integral~\cite{Bleistein1975Asymptotic,Chew1995waves,perez2017windowed,Brekhovskikh2013Acoustics}
\begin{eqnarray}
H_j(\ner,\ner')&=&\frac{1}{4\pi } \int_{\mathit{SC}} \frac{p_j(\xi,\ner')}{q(\xi)}\e^{|\ner|\phi(\xi)}\de \xi.\label{eq:asymp_int_R}
 \end{eqnarray}
 Here, letting $\gamma_j(\xi)=\sqrt{\xi^2-k_j^2}$, $j=1,2,3$, we have
 set
\begin{subequations}\begin{eqnarray}
\phi(\xi) &= & i\xi\cos\theta-\gamma_1(\xi)\sin\theta,\label{eq:phase_function}\\
p_1(\xi,\ner') &=&\lf\{R_{12}(\xi)+R_{23}(\xi)\e^{-2\gamma_2(\xi)d_2}\rg\}\frac{\e^{-i\xi x'-\gamma_1(\xi) y'}}{ \gamma_1(\xi)}\\
p_2(\xi,\ner')&=&\lf\{1-R_{12}(\xi)\rg\}\lf\{1+R_{23}(\xi)\e^{-2\gamma_2(\xi)(d_2+y')}\rg\}\frac{\e^{-i\xi x'+\gamma_2(\xi)y'}}{ \gamma_2(\xi)}\\
p_3(\xi,\ner')&=&\lf\{1-R_{12}(\xi)\rg\}\lf\{1-R_{23}(\xi)\rg\}\e^{-\gamma_2(\xi)d_2}
\frac{\e^{-i\xi x'+\gamma_3(\xi)(y'+d_2)}}{\gamma_3(\xi)}\\
q(\xi) &= &1+R_{12}(\xi)R_{23}(\xi)\e^{-2\gamma_2(\xi)d_2},
\end{eqnarray}
 and \begin{equation}
R_{ij} (\xi)=\frac{\gamma_i(\xi)-\nu_i\gamma_j(\xi)}{\gamma_i(\xi)+\nu_i\gamma_j(\xi)},\quad i,j=1,2,3.
\end{equation} \label{eq:asymp_int_2}\end{subequations} 

\begin{figure}[!ht]
  \centering
      \includegraphics[width=0.9\textwidth]{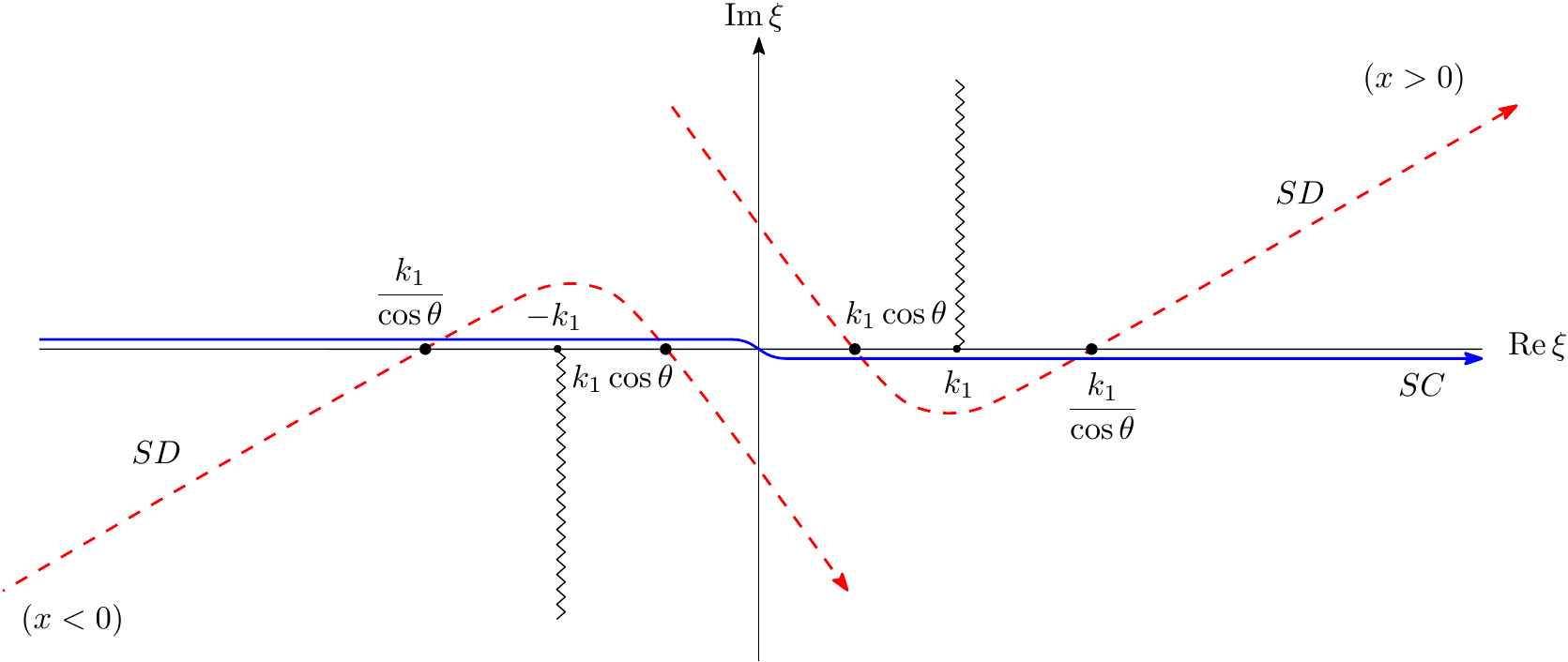}
      \caption{Sommerfeld contour $\mathit{SC}$ used
        in~\eqref{eq:asymp_int_R} (solid blue curve) and related
        steepest descent path $\mathit{SD}$ (dashed red curve).}
\label{fig:steepest_descent_reflected}
\end{figure}

The determination of physically admissible branches of the functions
$\gamma_j(\xi) =\sqrt{\xi^2-k_j^2}= \sqrt{\xi-k_j}\sqrt{\xi+k_j}$
requires adequate selection of branch cuts. Relevant branches, which
must be selected to insure the Green function satisfies outgoing
radiation condition for the layered structure, are given by $-3\pi/2
\leq \arg(\xi-k_j)<\pi/2$ for $\sqrt{\xi-k_j}$ and
$-\pi/2\leq\arg(\xi+k_j)<3\pi/2$ for $\sqrt{\xi+k_j}$. The branch cut
stemming from the point $\xi=k_1=k_3$ is in fact the only branch cut
in the domain of definition of the functions $p_j(\xi,\ner')/f(\xi)$
($j=1,2,3$), as it can be shown that these are even functions of
$\gamma_2$~\cite{Chew1995waves}. The branch cuts and Sommerfeld
contour~${\mathit{SC}}$ utilized in~\eqref{eq:asymp_int_R} are depicted in
Figure~\ref{fig:steepest_descent_reflected}.

As suggested above, in order to obtain the far-field form of the layer
Green function $H^f$ we resort to the method of steepest
descents~\cite{Bleistein1975Asymptotic}. Analysis of the phase
function $\phi$~\eqref{eq:phase_function} readily shows that there is
only one saddle point on the real axis at $\xi_0=k_1\cos\theta$ and
that the path of steepest descent $\mathit{SD}$ that passes through that point,
which is given by the expression $\imag \phi(\xi)=k_1$, also
intersects the real axis at $\xi=k_1/\cos\theta$. Furthermore, from
the definition of the function $\gamma_1$ it can shown that
$$
\imag \xi = \frac{|\cos\theta|}{\sin\theta}\,\real \xi -\frac{k_1}{\sin\theta}\quad\mbox{as}\quad |\xi|\to\infty,
$$
for $\xi$ on $\mathit{SD}$. This information suffices to sketch the paths of steepest descent that are displayed in~Figure~\ref{fig:steepest_descent_reflected}. 

In order to produce asymptotic expansions of the
integrals~\eqref{eq:asymp_int_R} we then proceed to deform the
Sommerfeld contour $\mathit{SC}$ to the steepest descent contour
$\mathit{SD}$ (Figure~\ref{fig:steepest_descent_reflected}). Considering the saddle point at $\xi_0$
and taking into account the poles of the integrand
$p_j(\xi,\ner')/q(\xi)$ at the points $\xi_p$ which are enclosed by
the curves $\mathit{SD}$ and $\mathit{SC}$, we obtain the following
expression for the far-field form of the Green function for the
two-dimensional slab:
\begin{equation}\label{hf_final}\begin{split}
  H^f(\ner,\ner' ) =&~\frac{i}{2}\sum_{\xi_p\in I}\underset{\xi=\xi_p}{\operatorname{Res}}\lf(\frac{p_j(\xi,\ner')}{q(\xi)}\e^{|\ner|\phi(\xi)}\rg)+\\
  &~\frac{1}{4\pi}\frac{p_j(\xi_0,\ner')}{q(\xi_0)}\sqrt{\frac{2\pi}{|\ner||\phi''(\xi_0)|}}\e^{|\ner|\phi(\xi_0)-i\pi/4} + O(|\ner|^{-3/2}),\quad (\ner'\in D_j)
\end{split}\end{equation}
as $|\ner|\to\infty$. Note that for $\cos\theta>0$
(resp. $\cos\theta<0$) only the real poles contained in the set
$I=(0,k_1\cos\theta)\cup(k_1/\cos\theta,\infty)$ (resp.
$I=(-\infty,k_1/\cos\theta)\cup(k_1\cos\theta,0)$) produce
contributions which do not decay exponentially.

Clearly, as indicated above, the far field asymptotics $H^f$ of the
layer Green function $H$, and thus its normal derivative, can be
evaluated inexpensively by means of a simple explicit expressions.

\section{Numerical examples}\label{sec:numerical_examples_ML}
This section presents a set of two- and three-dimensional numerical
examples that demonstrate the character of the proposed multi-layer
WGF methodology. For the sake of definiteness a window
function~$w_A$~\eqref{window_fnct} with $c=0.7$ was used in all
cases. Numerical errors were evaluated by resorting to
numerical-convergence studies and/or increases in the window-size
$A$. As additional references, in some cases adequately accurate
solutions obtained by the Sommerfeld layer-Green-function (LGF)
method~\cite{PerezArancibia:2014fg,perez2017windowed} (with accuracy evaluated by means
of convergence studies) were used to evaluate the accuracy of the WGF
approach. Brief indications will be provided when necessary to
indicate which method is being used in each case. The two-dimensional
results were obtained via solution of the integral equation
system~\eqref{eq:transmission_system_window} by means of the Nystr\"om
method described in \cite[Section 3.5]{COLTON:2012}. The
three-dimensional solutions, in turn, were obtained by means of the
algorithm presented in~\cite{Bruno:2001ima}.


\begin{figure}[ht!]
\centering
\includegraphics[width=0.4\textwidth]{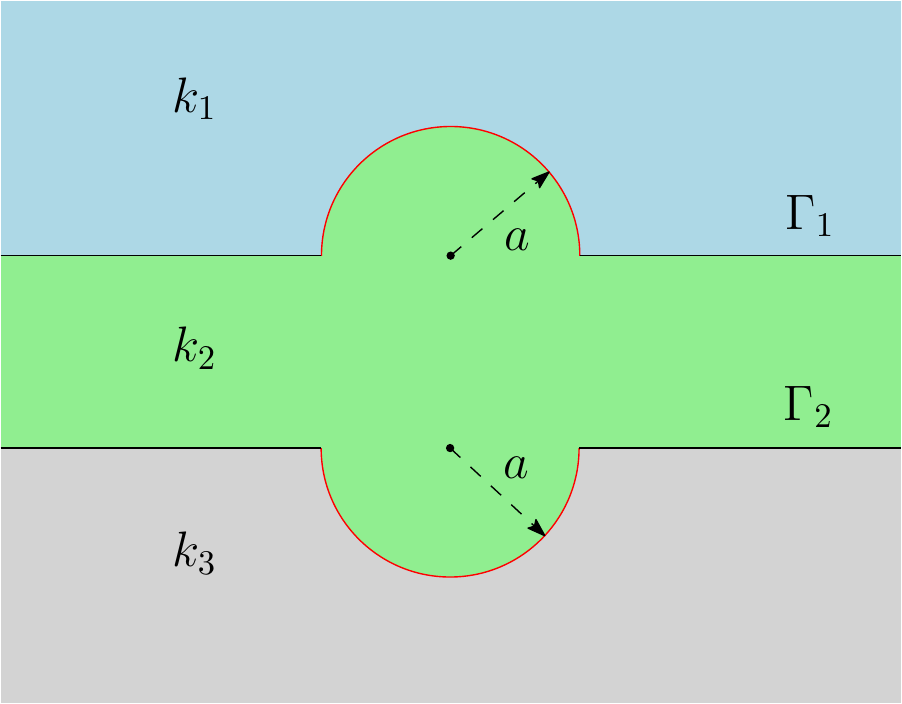}\qquad
\includegraphics[scale=0.53]{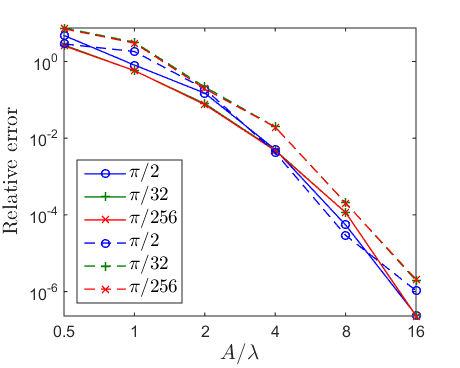}
\caption{Left. Structure utilized in the numerical examples presented
  in Section~\ref{sec:numerical_examples_ML}. Right. Relative errors
  in log-log scale in the integral densities resulting from numerical
  solution of~\eqref{eq:transmission_system_window} for the structure
  depicted on the left panel, by means of the WGF method, for various
  window sizes and angles of incidence---including extremely shallow
  incidences.  The WGF method computes integral densities with
  super-algebraically high accuracy uniformly for all incidences.}
  \label{fig:convergence_3L}
\end{figure}

Our first example concerns the structure depicted in the left
portion of Figure~\ref{fig:convergence_3L}, in which semi-circular
defects of radii $a=1$ are placed at the planar interfaces
$P_1=\R\times\{0\}$ and $P_2=\R\times\{-3/2\}$ of a three-layer medium
with wavenumbers $k_1=10$, $k_2=20$ and $k_3=30$.  The right portion
of Figure~\ref{fig:convergence_3L} displays the maximum relative
errors (in log-log scale) in the total field produced by the WGF
method on the surface of the semi-circular defects (the curves marked
in red in Figure~\ref{fig:convergence_3L} left) for various windows
sizes $A>0$ and incidences~$\alpha$.  The number of quadrature points
was selected in such a way that for any given $A>0$ the Nystr\"om
discretization error in the integral equation solution is not larger
than $10^{-9}$. The WGF solution obtained for $A=32\lambda$ is
utilized as the reference for the error estimation.  As it can be
inferred from the error curves displayed in
Figure~\ref{fig:convergence_3L}, super-algebraic convergence is
observed as $A$ increases. In particular, these results demonstrate
the uniformly fast convergence exhibited by the WGF method as the
incidence angles approach grazing.

\begin{table}[h!]
\centering
\begin{tabular}{c|c|c|c|c|c|c|c|c|c|c}
    \toprule
&\multicolumn{5}{c|}{WGF method} &\multicolumn{5}{c}{LGF method} \\
\hline     $\kappa$ &$2$&$4$&$8$&16&32 &$2$&$4$&$8$&16 &32\\
\hline    \scalebox{0.9}{Number of} \vspace{-0.1cm}&  \multirow{2}{*}{1232} &\multirow{2}{*}{1272}&\multirow{2}{*}{1348}& \multirow{2}{*}{1496}&\multirow{2}{*}{1800}&\multirow{2}{*}{68} &\multirow{2}{*}{148}&\multirow{2}{*}{300}&\multirow{2}{*}{596}&\multirow{2}{*}{1204}\\ \scalebox{0.9}{unknowns}&&&&&&&&&\\
\hline       \scalebox{0.9}{Matrix}\vspace{-0.1cm} &\multirow{2}{*}{3.44}  &\multirow{2}{*}{3.53} &\multirow{2}{*}{3.98} &\multirow{2}{*}{5.78}&\multirow{2}{*}{7.29} &\multirow{2}{*}{6.49} & \multirow{2}{*}{22.15} & \multirow{2}{*}{82.86}&\multirow{2}{*}{319.46}&\multirow{2}{*}{$1900$}\\ \scalebox{0.9}{construction (s)}&&&&&&&&&\\
\bottomrule
\end{tabular}
\caption{\label{tb:WGF_multi_times} Computing times required by the WGF and  LGF methods to construct the system matrices  for the numerical solution of  the problem of scattering of a plane-wave by a semi-circular cavity or radius $a=1$ on a three-layer medium with wavenumbers $k_1=\kappa$, $k_2=2\kappa$ and $k_3=3\kappa$, with $\kappa=2^j$, $j=1,\ldots 5$. All the two-dimensional 
  runs reported in this paper were performed using a Matlab implementation of our algorithms in a MacBook Air laptop (early 2014
  model).}
\end{table}

In order to compare the computational cost of the LGF
method~\cite{PerezArancibia:2014fg} and proposed WGF method for a
given accuracy, we consider a planar three-layer structure similar to
those considered previously, but now containing only one surface
defect: a semi-circular cavity of radius $a=1$ at~$P_1$. (The use of a
single defect reduces somewhat the LGF cost which seemed inordinately
large for the two-defect problem.) A plane-wave $u^\inc$ with $\alpha
=-\pi/6$ illuminates the structure. Five sets of wavenumbers given by
$k_1=\kappa,$ $k_2=2\kappa$ and $k_3=3\kappa$ with $\kappa=2^j$,
$j=1,\ldots,5$ are considered. The resulting problems of scattering
are then solved by employing a Nystr\"om discretization of the WGF
equations~\eqref{eq:transmission_system_window}, and a numerical
version of~\eqref{eq:simpler_near_field_evaluation_NL} is used to
evaluate near-fields. The same problem of scattering is then solved,
with a relative error not larger than $10^{-4}$, by means of a
generalization to the present three-layer case, of the two-layer LGF
method presented in~\cite{PerezArancibia:2014fg} (see
also~\cite{perez2017windowed}). The reference solution used to
estimate the accuracy of the LGF solution is obtained by solving the
resulting LGF integral equation with an error not larger $10^{-9}$
(this accuracy is achieved by utilizing a large number of Nystr\"om
quadrature points and evaluating the layer Green function with an
error not larger than $10^{-10}$).

\begin{figure}[ht!]
  \centering
	\includegraphics[width=0.4\textwidth]{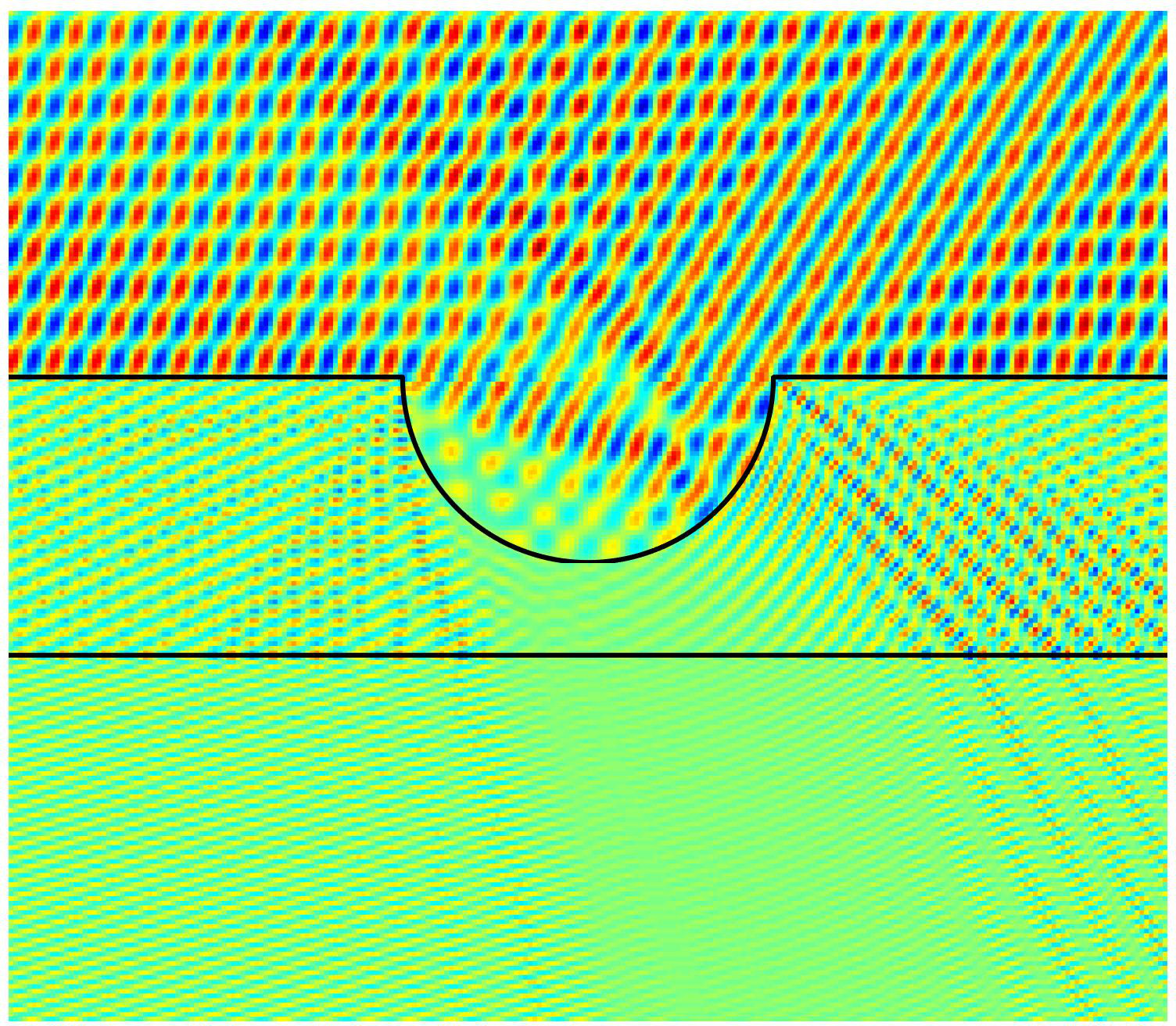}\qquad
	\includegraphics[width=0.4\textwidth]{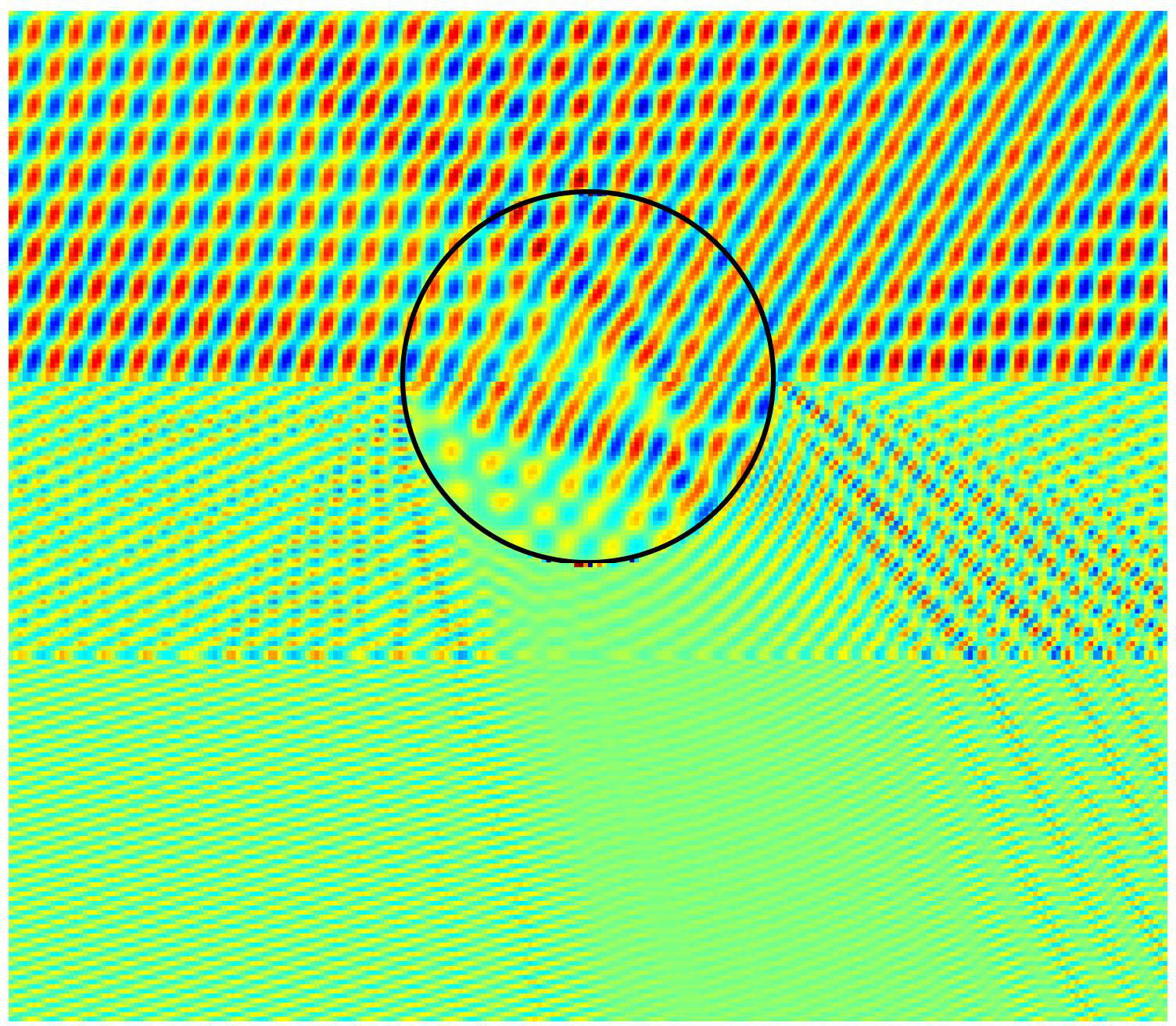}
        \caption{Real part of the total near fields obtained (with
          errors of the order of $10^{-4}$) by means of the WGF (left)
          and LGF (right) methods for the problem of scattering of a
          plane wave by a semi-circular cavity of radius $a=1$ on the
          top interface of a three-layer medium with wavenumbers
          $k_1=32$, $k_2=64$ and $k_3=96$ and incidence angle
          $\alpha=-\pi/6$. The respective integral-equation curves are
          shown in black. The WGF solution with $A=16\lambda$
          (resp. the LGF solution) was produced in total computing
          time of 62 secs (resp. $7.8\cdot 10^4$ secs).}
  \label{fig:near_field_3L_comparison}
\end{figure}

Table~\ref{tb:WGF_multi_times} displays the computing times needed by
both methods to construct the system matrices. In order to allow for a
fair comparison of the computing times and the field values on the
surface defect, the same set of quadrature points is utilized to
discretize the currents on the surface of the cavity in each case. The
number of quadrature points was increased in direct proportion to the
value of $\kappa$. The maximum of the absolute value of the difference
between the LGF and WGF solutions (using $A=8\lambda$) on the surface
of the defect is no larger than $10^{-4}$ in all the examples
considered. Remarkably, in the $\kappa=32$ case the proposed WGF
method is 260 times faster than the LGF method.

Figure~\ref{fig:near_field_3L_comparison} presents a comparison of the
near fields obtained by means of the WGF and LGF methods for some of
the test cases considered in Table~\ref{tb:WGF_multi_times}. The first
and second columns in Figure~\ref{fig:near_field_3L_comparison}
display the real-part of the total near-fields produced by the WGF
method (1st column) and by the LGF method (2nd column) respectively
for $\kappa=32$. The fields are evaluated in the rectangular region
$[-3,3]\times[-7/2,2]$ at an uniform grid of $280\times 200$
points. Note that, as it follows from consideration of the figure
captions, the WGF near field evaluation procedure is up to 1200 times
faster than the corresponding LGF near field evaluation procedure---in
spite of the fact that a (larger) window size $A=16\lambda$ had to be
used to produce accurate near fields throughout the plotted region.

\begin{figure}[ht!]
  \centering
   \subfloat[$\alpha = -\pi/2$.
                \label{fig:pi2_ff}]{ \includegraphics[width=0.48\textwidth]{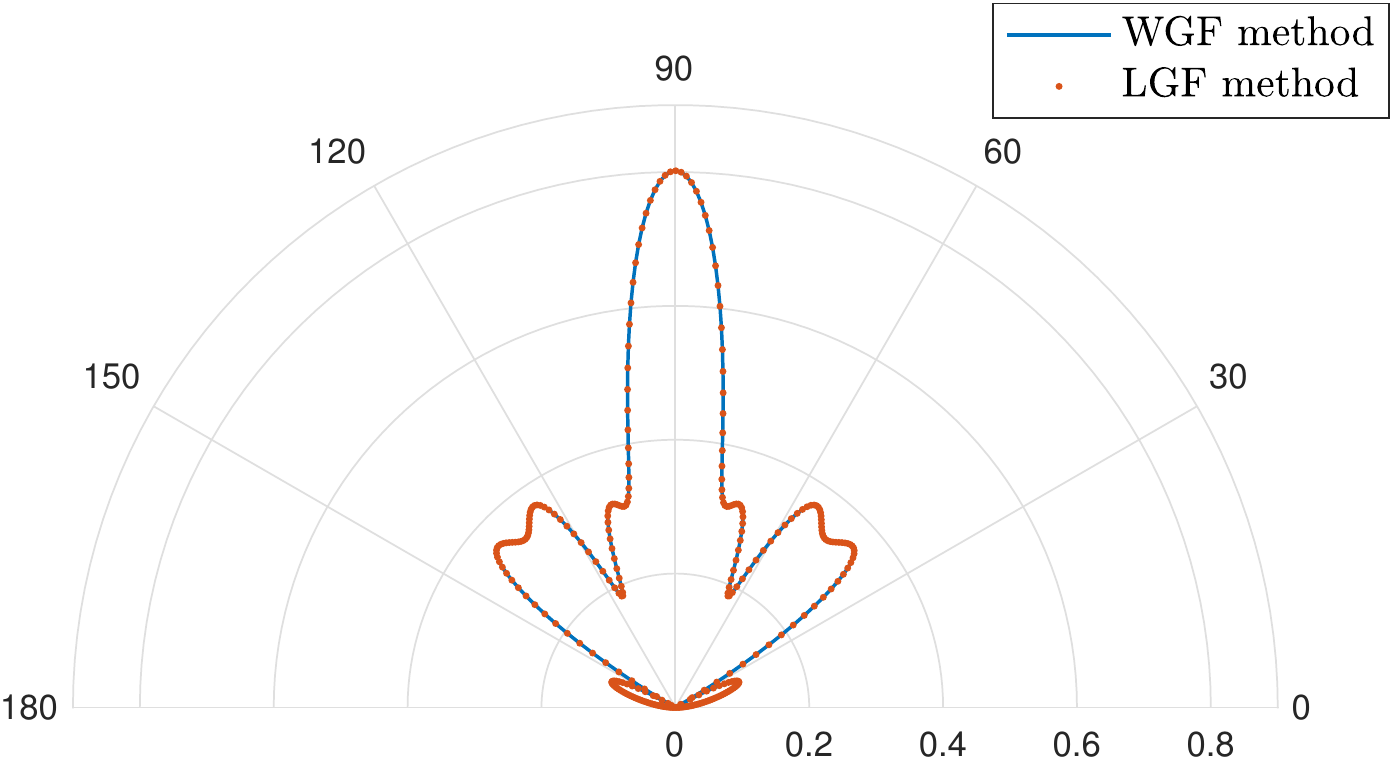}}                                         
 \subfloat[$\alpha = -\pi/6$.
                \label{fig:pi6_ff}]{ \includegraphics[width=0.48\textwidth]{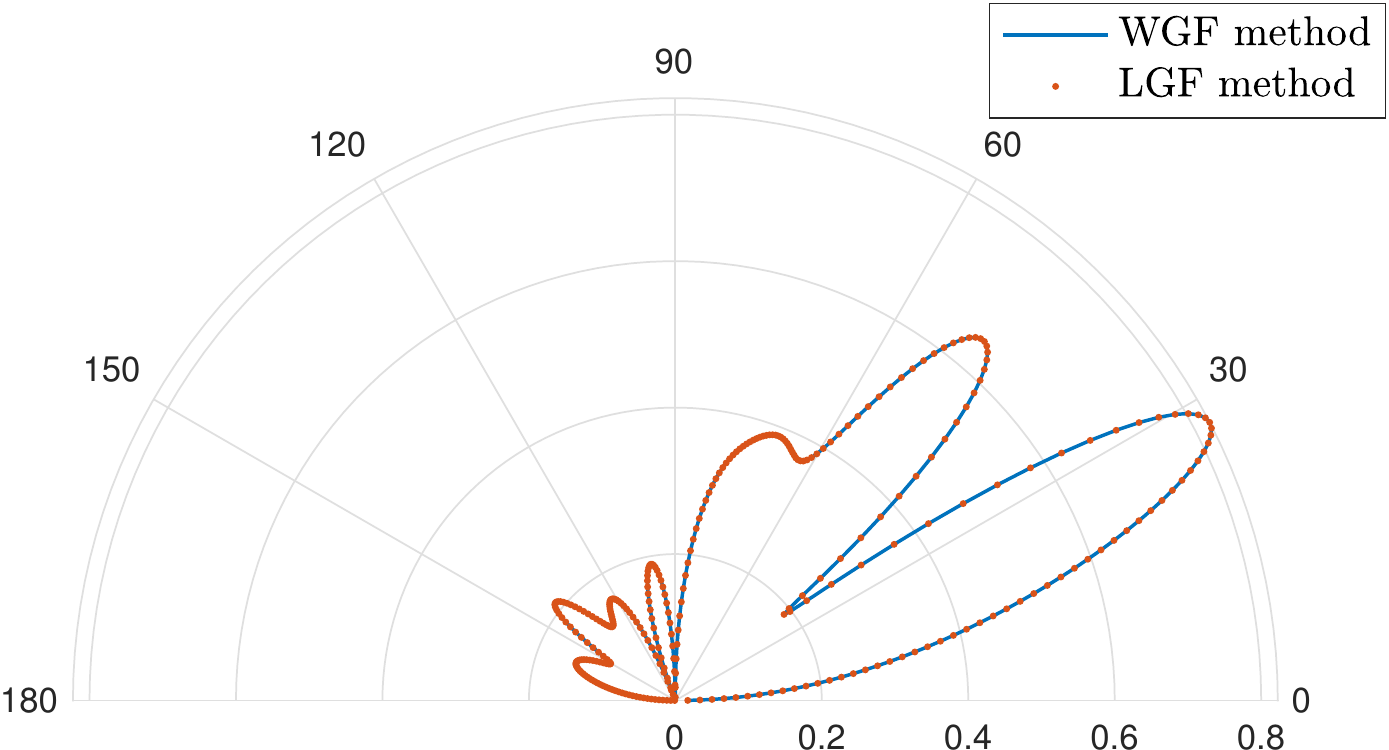}}
                \caption{WGF (blue solid line) and LGF~(red dotted
                  line) far-field patterns obtained for the problems
                  of scattering of a semi-circular cavity in a slab
                  with wavenumbers $k_1=k_3=10$ and $k_2=15$ for two
                  different incidence angles.}
  \label{fig:ff_field_3L}
\end{figure}

Figure~\ref{fig:ff_field_3L}, in turn, compares the far-field patterns
\begin{equation}\label{u_infty}
  u_{\infty}(\hat\ner) =\lim_{|\ner|\to\infty} \sqrt{|\ner|}\e^{-ik_1|\ner|}u(|\ner|\hat\ner),\qquad \hat\ner=\frac{\ner}{|\ner|}=(\cos\theta,\sin\theta),\quad \theta\in(0,\pi),
\end{equation}
produced by the WGF and LGF algorithms for a semi-circular cavity in a
three-layer medium with wavenumbers $k_1=k_3=10$ and $k_2=15$. The WGF
far-field pattern (blue solid line) was obtained by letting $u=u^f$
in~\eqref{u_infty}, where $u^f$ is given by~\eqref{eq:ext_field_ff}
with WGF defect fields $u^d = u^{d,w}_j$ in $\Omega_j$, $j=1,\ldots,N$
(equation~\eqref{eq:simpler_near_field_evaluation_NL}). The
corresponding LGF far-field pattern (red dots) was obtained on the
basis of a highly accurate LGF solution together with the far-field
asymptotics of the layer Green function~\cite{PerezArancibia:2014fg}.
We have verified that, as expected, the accuracy of the WGF far-field
patterns is comparable to the accuracy of the corresponding defect
fields $u^{d,w}_j$ within the region~$\tilde \Omega_A$.
 
 
\begin{figure}[ht!]
        \centering
 \includegraphics[height=0.35\textwidth]{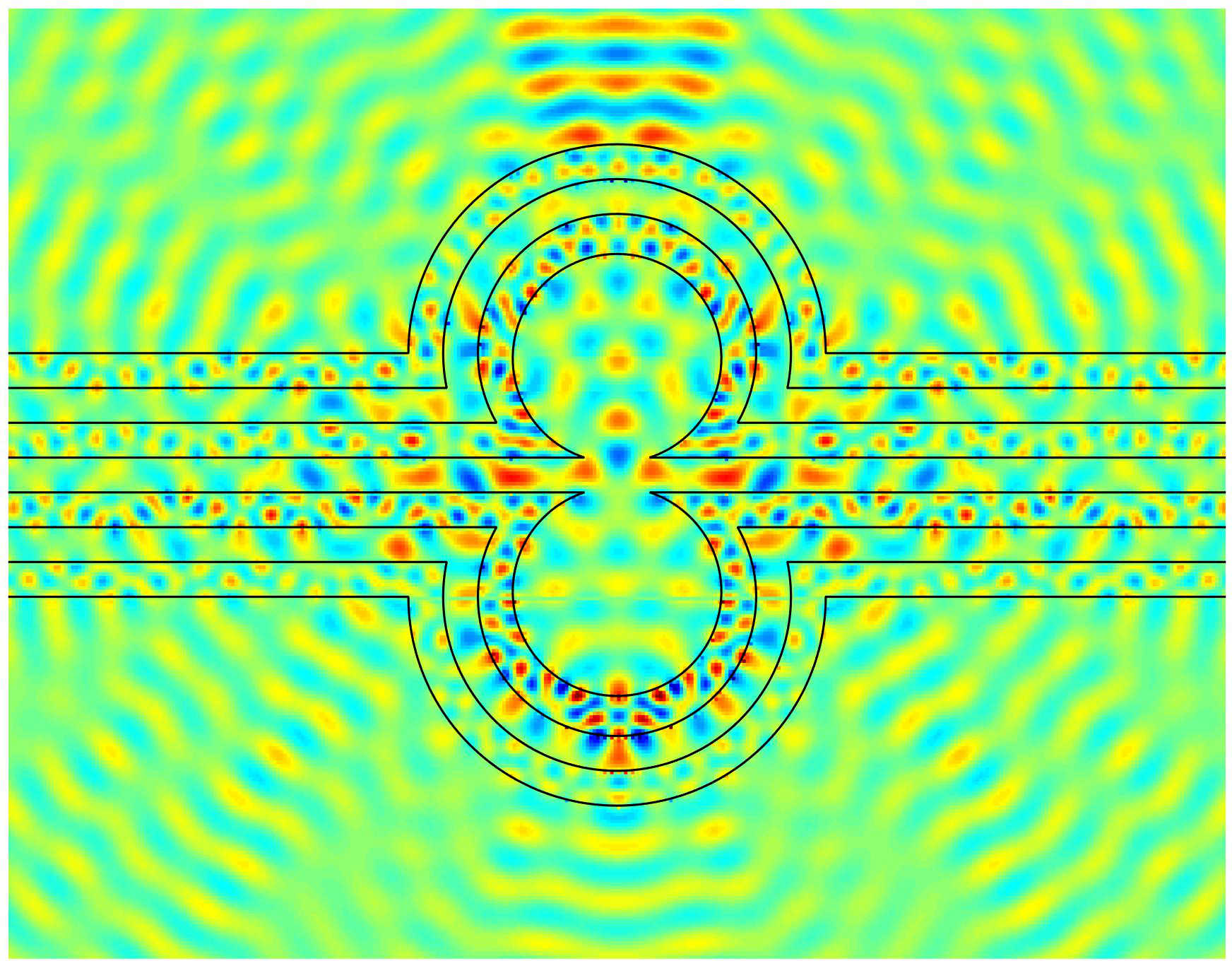}\quad \includegraphics[height=0.35\textwidth]{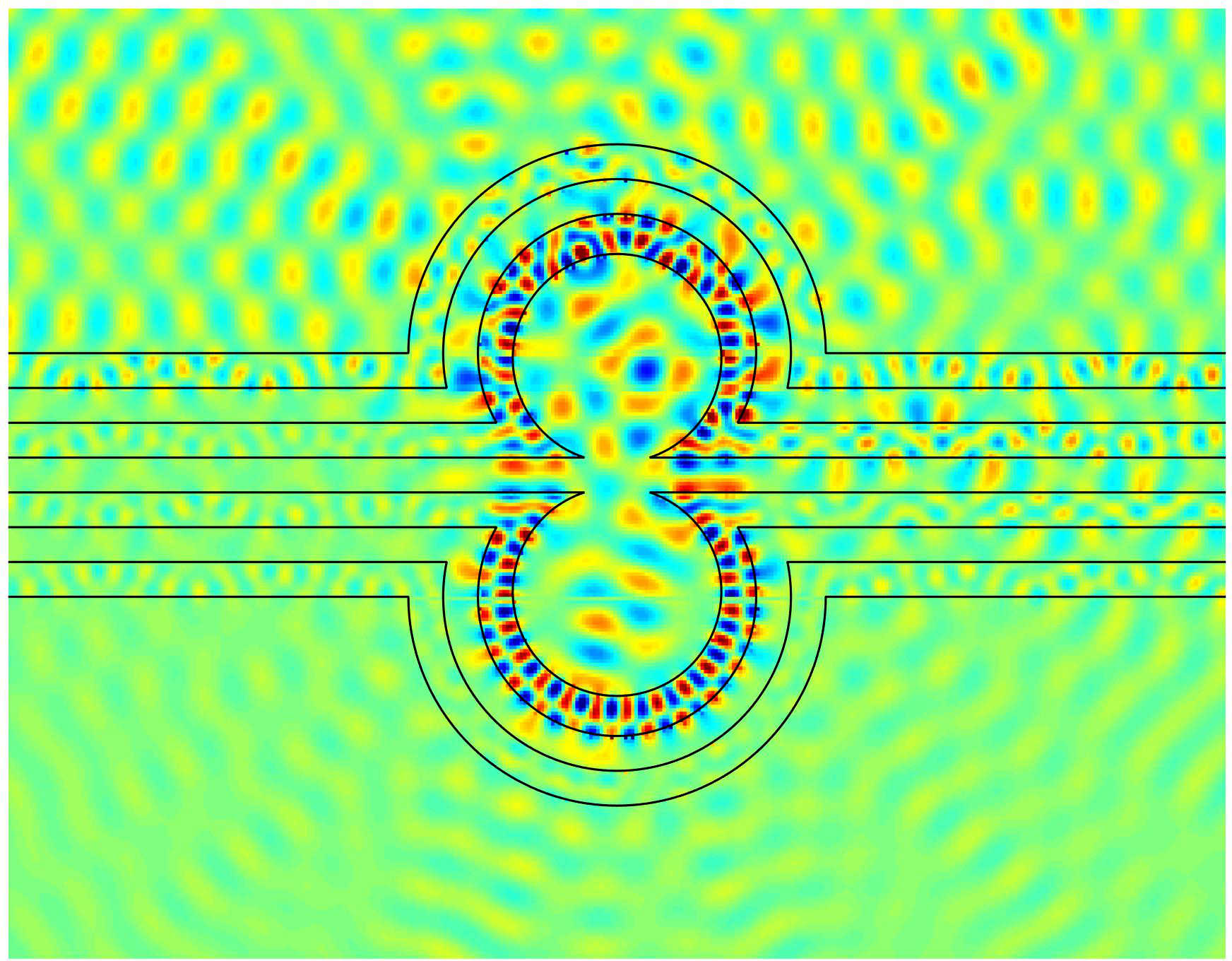}
 \caption{Real part of the total field for the problem of scattering
   of a plane wave impinging on a layered medium composed by 9 layers:
   $k_{2j-1}=15$, $j=1,\ldots,5$ and $k_{2j}=30$, $j=1,\ldots,4$ and
   $P_j=\R\times\{(j-1)/5\}$,
   $j=1,\ldots,8$. }\label{fig:corner_9_layers}
\end{figure}

Figure~\ref{fig:corner_9_layers} displays near fields resulting from
the WGF method, with window size $A=12\lambda$, for a structure
consisting of nested circular surface defects in a nine-layer medium
with planar interfaces $P_j=\R\times\{(j-1)/5\}$, $j=1,\ldots,8$. The
corresponding wavenumbers are $k_{2j-1}=15$ for $j=1,\ldots, 5$ and
$k_{2j}=30$ for $j=1,\ldots 4$. The structure is illuminated by
plane-waves with two different incidence angles.  A 112-second overall
computing time sufficed to evaluate each one of the two near fields
displayed. Note the resonance that takes place in the third upper and
lower rings in Figure~\ref{fig:corner_9_layers} right.

\begin{figure}[ht!]
        \centering
  \includegraphics[width=0.9\textwidth]{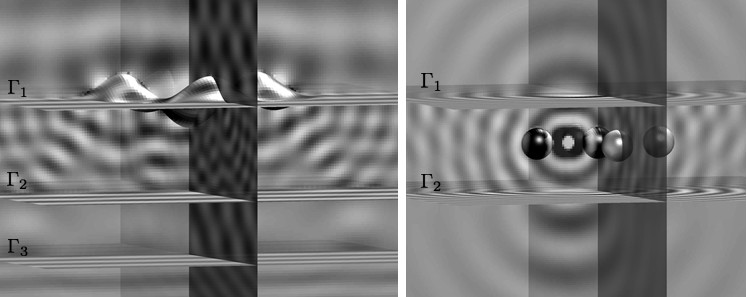}
  \caption{Three-dimensional total fields (real parts shown) produced
    by the WGF approach. Left. Scattering of a plane-wave by a surface
    defect in a four-layer medium with wavenumbers
    $k_1=k_3=4$ and $k_2=k_4=8$. Right. Scattering of a point-source
    field by an array of five spheres in a three-layer medium with
    wavenumbers $k_1=k_3=4$ and $k_2=8$. The absolute errors in the
    surface fields displayed are no larger than $4\cdot 10^{-4}$ for
    corresponding maximum fields of order one.}\label{fig:3d}
\end{figure}

Figure~\ref{fig:3d}, finally, presents applications of the WGF
methodology to the problem of scattering by three-dimensional
structures in presence of layer media. The two-dimensional
descriptions presented in Sections~\ref{sec:ml_prelim}
through~\ref{sec:far_field_complicated_NL} extend directly to the
present three-dimensional context.

\section*{Acknowledgements} 
This work was supported by NSF and AFOSR through contracts DMS-1411876
and FA9550-15-1-0043, and by the NSSEFF Vannevar Bush Fellowship under
contract number N00014-16-1-2808.

\appendix


\section{Appendix: Integral representation based on non-windowed
  free-space Green functions}\label{eq:green_multi_layer}

This section presents an integral representation formula, based on the
{\em free-space Green function}, for fields of the form
$v(\ner)=v_j(\ner)$ for $\ner\in\Omega_j$, $j=1,\ldots,N$, where,
letting $u^d_j$ and $\tilde u^p_j$ be defined in~\eqref{eq:def_field}
and~\eqref{eq:multi_layer_plane_wave_2}, respectively, we have either
\begin{equation}\label{eq:decomp_to_field}
  v_j = u^d_j + \tilde u^p_j, \quad v_j=u^d_j\quad\mbox{or}\quad v_j=\tilde u_j^p\quad \mbox{in}\quad \Omega_j.
\end{equation} 
The presentation is restricted to two-dimensional configurations. A
related (modified) representation, which can similarly be utilized for
all purposes necessary in this paper, can be obtained
analogously---albeit with certain additional considerations, as
detailed in~\cite{DeSanto:1998tq} in the three-dimensional sound-hard
case; cf. also~\cite{DeSanto:1997es}.  For simplicity, the
presentation is further restricted to three-layer structures, but the
extension to $N$-layer structures is straightforward.

\begin{figure}[h!]
\centering	
\includegraphics[scale=0.8]{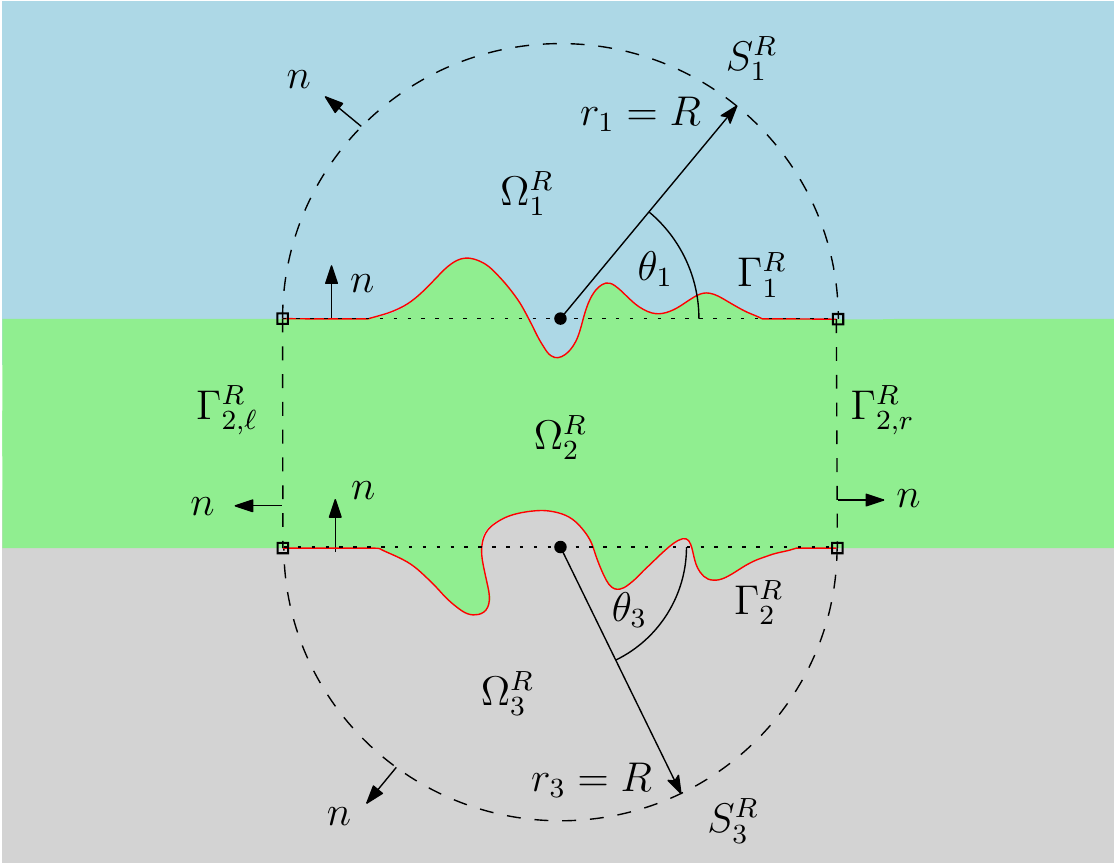}
\caption{Depiction of the various domains, boundaries and variables
  involved in the derivation of the integral representation
  formula~\eqref{eq:green_representation}. }\label{fig:green_domain}
\end{figure}

Our derivations utilize three local polar-coordinate systems, each one
of which is associated with one of the layers~$\Omega_j$. These
coordinate systems are centered at $(0,-d_1)$, $(0,-(d_{1}+d_{2})/2)$
and $(0,-d_{2})$ and, thus, the radial variables are given by
\begin{equation}
      r_1=\sqrt{x^2+(y+d_1)^2},\quad
      r_2=\sqrt{x^2+(y+(d_{1}+d_{2})/2)^2},\quad\mbox{and}\quad
      r_3=\sqrt{x^2+(y+d_{2})^2},
\label{eq:radius_multi}
\end{equation}
in terms of the global Cartesian coordinates $x$ and $y$, as
illustrated in Figure~\ref{fig:green_domain}. 

Additionally, some of the subsequent derivations utilize the
decomposition
\begin{equation}
\tilde u^p_j=u_j^{\uparrow}+u_j^{\downarrow} \quad\mbox{in}\quad \Omega_j,
\label{eq:total_flat_decomp}
\end{equation}
where letting $k_{1x} = k_1\cos\alpha$ and $k_{jy} =
\sqrt{k_j^2-k^2_{1x}}$, the up-going and down-going plane-waves
$u_j^{\uparrow}$ and $u_j^{\downarrow}$ are given by
\begin{equation}
\begin{split}
u_j^{\uparrow}(\ner) = p_j\e^{ik_{1x}x+ik_{jy}y}\quad \mbox{and} \quad
u_j^{\downarrow}(\ner) = q_j\e^{ik_{1x}x-ik_{jy}y},
\end{split}\label{eq:wave_up_down}
\end{equation}
respectively. Here the constants $p_j=\e^{2ik_{jy}d_j}A_j\widetilde
R_{j,j+1}$ and $q_j=A_j$ are expressed in terms of the amplitudes
$A_j$ and the generalized reflection coefficients $\widetilde
R_{j,j+1}$ defined in~\eqref{eq:ampl} and~\eqref{eq:gen_ref}.  Note
that $u_1^{\downarrow} = u^\inc$ and $u_3^{\uparrow}=0$.  The defect
field $v^d_j$, on the other hand, is given by~\cite{JerezHanckes:2012hk,Nosich:1994vr}
\begin{equation}
  v^d_j=\left\{\begin{array}{ccl}
      v^{\mathrm{rad}}_j+v^{\mathrm{gui}}_j &\mbox{in}&\Omega_j ,\quad j=1,3,\\
      v^{\mathrm{gui}}_2&\mbox{in}&\Omega_2,\quad j=2,
\end{array}\right.\label{eq:scatt_decomp}
\end{equation}
in terms of radiative and guided wave fields $v_j^{\mathrm{rad}}$ and
$v_j^{\mathrm{gui}}$ which, letting $\beta_2=1$ and $\beta_j=2/3$ for
$j=1,3$ (see~\eqref{eq:radius_multi}), verify
\begin{equation}
\lim_{r_j\to\infty}\sqrt{r_j}\left(\frac{\p v^\mathrm{rad}_j}{\p r_j} -ik_j v^\mathrm{rad}_j\right) = 0,\quad \quad j=1,3,\label{eq:Somm_layer}\end{equation}
and 
\begin{equation}\left\{\begin{array}{rcl}
      v^{\mathrm{gui}}_j(\ner) =\displaystyle \sum_{m=1}^{M_j}\alpha^m_{j} v^m_{j}(\ner)+ O\left(r_j^{-\beta_j}\right),&\\
      \displaystyle\left|\frac{\p v^\gui_j}{\p r_j}-i \sum_{m=1}^{M_j} \alpha^m_{j} \xi^m_{j} v^m_{j}\right| =\mathcal  O\left(r_j^{-\beta_j}\right)&\end{array}\mbox{ as } r_j \to\infty,\ (j=1,2,3).\right.\label{eq:rad_cond_modes}\end{equation}
Here  $v^m_j$ denote the guided modes
\begin{equation}
v_j^m(\ner) =\left\{\begin{array}{cl}
\left\{a_2^m \cosh\left(\gamma_2^m y\right)+b_2^m\sinh\left(\gamma_2^m y\right)\right\}\e^{i|x|\xi_2^m},& j=2 ,\medskip\\
 \e^{-\gamma_j^m |y|}\e^{i|x|\xi_j^m},& j=1,3,
\end{array}\right.
\label{eq:modes}\end{equation}
which are expressed in terms of the so-called propagation constants
$\xi_j^m>0$, and $\gamma_j^m = \sqrt{(\xi_j^m)^2-k_j^2}$, $m=1,\ldots,
M_j$.  The propagation constants $\xi_j^m$ equal the real poles
(sometimes called surface wave poles~\cite{Cho:2012hv,Cai:2000bl}) of
the corresponding three-layer Green function in spectral form. The
condition for the existence of the propagative modes in the inner
layer~$\Omega_2$ is $k_1<\xi_2^m<k_2$. For the outer layer $\Omega_1$
(resp. $\Omega_3$), on the other hand, we have $\xi_1^m=\xi_2^m$
(resp.~$\xi_3^m = \xi_{2}^m$) and the guided-mode condition is
$\xi_1^m>k_1$ (resp. $\xi_3^m>k_3$). Thus $v_1^m$ (resp. $v_3^m$)
corresponds to a surface wave that travels along the interface
$\Gamma_1$ (resp. $\Gamma_{2}$) and decays exponentially fast towards
the interior of $\Omega_1$ (resp. $\Omega_3$).

We are now in a position to derive the desired integral representation
for the total fields $v_j$ in~\eqref{eq:decomp_to_field}.  Our
derivations consider at first the bounded domains
\begin{equation}
B_R =  \left((-R,R)\times (-d_{2},-d_1)\right)\cup 
\left\{(x,y): r_1<R\right\}\cup\left\{(x,y): r_3<R\right\}\label{eq:domain_for_green}
\end{equation}
where $R>0$ is large enough that $B_R$ contains all of the surface
defects, as illustrated in Figure~\ref{fig:green_domain}.
Our bounded-domain calculations use the curves $\Gamma^R_{2,l}$,
$\Gamma^R_{2,r}$, $S^R_1$ and $S^R_3$ and corresponding normals $n$,
as depicted in Figure~\ref{fig:green_domain}. In order to facilitate
repeated use of Green's third identity in our derivations we
follow~\cite{DeSanto:1997es} and, letting
$G_{k_j}(\ner,\ner')=\frac{i}{4}H_0^{(1)}(k_j|\ner-\ner'|)$ we define,
for a given curve $C$,
\begin{equation}\label{Green-f}
  I_j\lf[v;C\rg](\ner) = \int_{C}\left\{\frac{\p G_{k_j}}{\p n_{\ner'}}(\ner,\ner')v(\ner')-G_{k_j}(\ner,\ner')\frac{\p v}{\p n}(\ner')\right\}\de s_{\ner'}.
\end{equation}
In what follows, finally, we make frequent use of the
$|\ner'|\rightarrow\infty$ asymptotic relations
\begin{equation}
\begin{array}{rcc}
G_k(\ner,\ner') &=&\displaystyle\frac{i\e^{-i\pi/4}}{\sqrt{8\pi k |\ner'|}}\e^{ik(|\ner'|- \ner\cdot\hat\ner')} \left\{1+O\left(|\ner'|^{-1}\right)\right\},\medskip\\
\displaystyle \nabla_{\ner'} G_k(\ner,\ner') &=&\displaystyle-\frac{\sqrt{k}\e^{-i\pi/4}}{\sqrt{8\pi  |\ner'|}}\e^{ik(|\ner'| - \ner\cdot\hat\ner')}\hat \ner'\left\{1+O\left(|\ner'|^{-1}\right)\right\}\\\end{array}\quad \left(\hat\ner'=\frac{\ner'}{|\ner'|}\right) \label{eq:asmyp_GF}
\end{equation} 
that follow directly from the corresponding asymptotic expressions for
the Hankel function~\cite{Lebedev:1965} together with easily verified
identity $|\ner-\ner'| = |\ner'|-\ner'\cdot\ner/|\ner'|
+O\lf(|\ner'|^{-1}\rg)$.


With reference to Figure~\ref{fig:green_domain}, and in view of
Green's third identity applied to $\Omega_1^R$ and its boundary
$\p\Omega_1^R = \Gamma_1^R\cup S_1^R$ we obtain the bounded-domain
integral representation
\begin{subequations}
\begin{equation}
 I_1\left[v_1;\Gamma_1^R\right](\ner) - I_1\left[v_1;S^R_1\right](\ner)=\left\{\begin{array}{cl} v_1(\ner), &\ner\in\Omega^R_1,\smallskip\\
 0,& \ner\in\R^2\setminus\Omega^R_1.
 \end{array}\right.\label{eq:green_1}
\end{equation}
Similarly, integrating over the domains $\Omega_2^R$ and $\Omega_3^R$,
whose boundaries are given by
$\p\Omega_2^R=\Gamma_2^R\cup\Gamma_{2,r}^R\cup\Gamma_1^R\cup\Gamma_{2,\ell}^R$
and $\p\Omega_3^R=\Gamma_2^R\cup S_3^R$, respectively, the
bounded-domain integral representations
\begin{equation}
 I_2\left[v_2;\Gamma_{2}^R\right](\ner)- I_2\left[v_2;\Gamma_{1}^R\right](\ner) - I_2\left[v_2;\Gamma^R_{2,\ell}\cup\Gamma^R_{2,r}\right](\ner) =\left\{\begin{array}{cl} v_2(\ner), &\ner\in\Omega^R_2,\smallskip\\
 0,& \ner\in\R^2\setminus\Omega^R_2,
 \end{array}\right.\label{eq:green_j}
\end{equation}
and 
\begin{equation}
 -I_3\left[v_3;\Gamma_{2}^R\right](\ner)  - I_3\left[v_3;S^R_3\right](\ner)=\left\{\begin{array}{cl} v_3(\ner), &\ner\in\Omega^R_3,\smallskip\\
 0,& \ner\in\R^2\setminus\Omega^R_3
 \end{array}\right.\label{eq:green_N}
\end{equation}\label{eq:green_truncated}\end{subequations}
result. In order to complete our calculations it suffices to evaluate
the limiting values as $R\to\infty$ for the various integral
quantities in~\eqref{eq:green_truncated} and for each one of the
functions $v_j$ in~\eqref{eq:decomp_to_field}. Since, in view of
equations~\eqref{eq:total_flat_decomp},~\eqref{eq:wave_up_down} and
\eqref{eq:scatt_decomp}, these functions can be expressed as linear
combinations of $u_j^\uparrow$, $u_j^\downarrow$, $v_j^{\mathrm{gui}}$
and $v_j^{\mathrm{rad}}$ in what follows we obtain the corresponding
limiting values for each one of these functions.  The desired
representation formulae~\eqref{eq:green_representation} as well as
their $N$-layer versions~\eqref{eq:green_representation_NL} then
follow directly from the limiting expressions thus found.


\noindent{\bf Case $j=2$ for $v_2 =u_2^\uparrow $,
  $v_2=u_2^\downarrow$, $v_2=v_2^{\mathrm{gui}}$ and
  $v_2=v_2^{\mathrm{rad}}$}. In view of the decay of the integral
kernels~\eqref{eq:asmyp_GF} and the fact that the total field $v_2$
remains bounded throughout $\Omega_2$ (as it follows from
equation~\eqref{eq:rad_cond_modes}), we conclude that the term $I_2$
involving the integral over $\Gamma^R_{2,r}\cup \Gamma_{2,\ell}^R$
tends to zero as $R\to\infty$. 


\noindent{\bf Case $j=1,3$ for $v_j=v_j^{\mathrm{rad}}$}.  In order to
estimate the terms $I_j$ that involve integrals over the semi-circular
curves $S_1^R$ and $S_3^R$, in turn, we note that for $\ner'\in S_j^R$
with $j=1$ and $j=3$ we have $|\ner'| = R+ O(1)$ and $\widehat{\ner'}=
(\cos\theta_j,\sin\theta_j)+ O(R^{-1})$ as $R\to\infty$---where the
the angles $\theta_j$ are as shown in
Figure~\ref{fig:green_domain}. Since $v^{\mathrm{rad}}_j$
in~\eqref{eq:scatt_decomp} for $j=1,3$ satisfies the Sommerfeld
radiation condition~\eqref{eq:Somm_layer}, utilizing standard
arguments~\cite{COLTON:1983} it can be shown that
$$
I_j\left[v_j^{\mathrm{rad}};S_j^R\right]=o\left(1\right),\quad j=1,3,\quad\mbox{as}\quad R\to\infty.
$$

\noindent{\bf Case $j=1,3$ for $v_j=v_j^{\mathrm{gui}}$}.  We now
consider the quantity $I_1\left[v^{\gui}_1;S_1^R\right]$, which,
according to~\eqref{eq:rad_cond_modes}, is given by a linear
combination of terms of the form $I_1\left[v_1^{m};S_1^R\right]$
($m=1,\ldots,M_j$), where $v_1^m(\ner) = \e^{-\gamma_1^m|y|+i\xi_1^m
  |x|}$ with $\gamma_1^m>0,$. From~\eqref{eq:asmyp_GF} and the fact
that $v_1^m(\ner') = \e^{-\gamma_1^m R\sin\theta_1+i\xi_1^m
  R|\cos\theta_1|}$ for $\ner'=R(\cos\theta_1,\sin\theta_1)\in S_1^R$,
$\theta_1\in[0,\pi]$, we obtain
$$
I_1\left[v_1^{m};S_1^R\right](\ner) \sim \sqrt{\frac{k_1R}{8\pi}}\e^{ik_1R-i\frac\pi4}\times\hspace{10cm}$$
$$\hspace{3cm}\int_{0}^{\pi}\left\{\frac{i\gamma_1^m\sin\theta_1+|\cos\theta_1|\xi_1^m}{k_1}-1\right\}\e^{-ik_1\ner\cdot\hat\ner'-R(\gamma_1^m\sin\theta_1-i\xi_1^m|\cos\theta_1|)}\de\theta_1
$$  as $R\to\infty$.  Therefore  
$$
\left| I_1\lf[v_1^{m};S_1^R\rg](\ner)\right| \leq  \sqrt{\frac{k_1R}{2\pi}}\left\{1+\frac{|\gamma_1^m|+\xi_1^m}{k_1}\right\}\int_{0}^{\pi/2}\e^{-\gamma_1^m R\sin\theta_1}\de\theta_1.
$$
The integral in the expression on the right-hand-side can be bounded
utilizing the inequality $\sin\theta\geq 2\theta/\pi$, $\theta\in
[0,\pi]$, and we thus conclude that $ I_1\lf[v_1^{m};S_1^R\rg]=
O\lf(R^{-1/2}\rg)$ as $R\to\infty$. Similarly, it can be shown that
$I_3\left[v_3^m;S_3^R\right] = O\left(R^{- 1/2}\right)$, and
consequently, in view of~\eqref{eq:rad_cond_modes}, we conclude that
$I_j\left[v_j^{\mathrm{gui}};S_j^R\right] = O\left(R^{-1/3}\right)$,
$j=1,3,$ as
$R\to\infty$.

\noindent {\bf Case $j=1$ for $v_1= u_1^\uparrow$ and $v_1=
  u_1^\downarrow$}. We consider the term
$I_1\lf[u_1^{\uparrow};S_1^R\rg]$ first, where $u^\uparrow_1(\ner) =
\e^{ik_{1x}x+ik_{1y}y}$, or, in polar coordinates,
$u^{\uparrow}_1(\ner') =\e^{ik_1R\cos(\theta_1+\alpha)}$ for $\ner'\in
S_1^R$. Then, integration by parts yields
\begin{align*}
I_1\lf[u_1^{\uparrow};S^R_1\rg](\ner) \sim&\ \sqrt{\frac{k_1R}{8\pi}}\e^{ik_1R-i\frac\pi4}\int_{0}^{\pi}\left\{\cos(\theta_1+\alpha)-1\right\}\e^{-ik_1  \ner\cdot\hat\ner'} \e^{ik_1R\cos(\theta_1+\alpha)}\de\theta_1\\
=&-\frac{\e^{ik_1R+i\frac\pi4}}{\sqrt{8\pi k_1 R}}\int_{0}^{\pi}\frac{\sin(\theta_1+\alpha)}{1+\cos(\theta_1+\alpha)}\e^{-ik_1 \ner\cdot\hat\ner'} \frac{\de }{\de \theta_1}\e^{ ik_1R\cos(\theta_1+\alpha)}\de\theta_1\\
=&-\frac{\e^{ik_1R+i\frac\pi4}}{\sqrt{8\pi k_1 R}}\left\{\left.\frac{\sin(\theta_1+\alpha)\e^{ik_1 (-\ner\cdot\hat\ner'+R\cos(\theta_1+\alpha))}}{1+\cos(\theta_1+\alpha)}\right|_{0}^{\pi}-\right.\\
&\left.\int_{0}^{\pi}\e^{ ik_1R\cos(\theta_1+\alpha)} \frac{\de }{\de \theta_1}\left(\frac{\sin(\theta_1+\alpha)\e^{-ik_1\ner\cdot\hat\ner'}}{1+\cos(\theta_1+\alpha)}\right)\de\theta_1\right\}=O\left(R^{-1/2}\right)
\end{align*}
as $R\to\infty$.  In order to deal with the term
$I_1\lf[u_1^{\downarrow};S^R_1\rg]$ with $u^\downarrow_1(\ner) =
\e^{ik_{1x}x-ik_{1y}y}$ we proceed similarly. Using the polar form
$u^{\downarrow}_1(\ner')= \e^{ik_1R\cos(\theta_1-\alpha)}$ for
$\ner'\in S_1^R$ of the down-going wave we obtain
\begin{align*}
&I_1\left[u_1^{\downarrow};S^R_1\right](\ner) \sim \sqrt{\frac{k_1R}{8\pi}}\e^{ik_1R-i\frac\pi4}\!\!\int\limits_{0}^{\pi}\left\{\cos(\theta_1-\alpha)-1\right\}\e^{-ik_1\ner\cdot\hat\ner'} \e^{ik_1R\cos(\theta_1-\alpha)}\de\theta_1.
\end{align*}
Note that since $\alpha\in (-\pi,0)$ we have
$0<\theta_1-\alpha<2\pi$. Thus, there is only one point of stationary
phase within the domain of integration at $\theta_1=\alpha+\pi$.  A
straightforward application of the method of stationary
phase~\cite{Bleistein1975Asymptotic} then yields
\begin{align*}
&I_1\left[u_1^{\downarrow};S^R_1\right](\ner)=-\e^{ik_1 (x\cos\alpha+y\sin\alpha)}+O\lf(R^{-1/2}\rg)=-u^\inc(\ner)+ O\lf(R^{-1/2}\rg)\quad\mbox{as}\quad R\to\infty.
\end{align*}
(Notice that integrating by parts yields that the limit points of the integral give rise to contributions that decay as $R^{-1}$.)

\noindent {\bf Case $j=3$ for $v_j= u_j^\downarrow$}.  This case
concerns the term $I_3\left[u^{\downarrow}_3;S_3^R\right]$ with
$u^\downarrow_3(\ner') =\e^{ik_{3x}x-ik_{3y}y}$, where
$k_{3x}=k_1\cos\alpha$ and $k_{3y}=\sqrt{k_3^2-k_{3x}^2}$. We
distinguish three possible cases, namely: (a)~$k_3<k_1|\cos\alpha|$,
(b)~$k_3=k_1|\cos\alpha|$ ($k_3=-k_1\cos\alpha$ for
$\alpha\in(-\pi,-\pi/2]$ or $k_3=k_1\cos\alpha$ for
$\alpha\in(-\pi/2,0)$), and (c)~$k_3>k_1|\cos\alpha|$.  Since in
case~(a) we have $k_{3y} = i\sqrt{k_1^2\cos^2\alpha-k_3^2}$, a
calculation completely analogous to the one carried in the estimation
of the term $I_1\left[u_1^{m};S_1^R\right]$ shows that
$I_3\left[u^{\downarrow}_3;S_3^R\right] = O\left(R^{- 1/2}\right)$.
In case~(b), in turn, we first consider $\alpha\in(-\pi/2,0)$. Under
this assumption $u_3^{\downarrow}(\ner') = \e^{ik_3R\cos\theta_3}$ for
$\ner'\in S^R_3$, and consequently
\begin{equation*}
\begin{split}
I_3\left[u_3^{\downarrow};S_3^R\right](\ner) \sim&\ \sqrt{\frac{k_3R}{8\pi}}\e^{ik_3R-i\frac\pi4}\int_{-\pi}^{0}\left\{\cos\theta_3-1\right\}\e^{-ik_2(d_{2}\sin\theta_3+\ner\cdot\hat\ner')+iRk_3\cos\theta_3}\de\theta_3.
\end{split}
\end{equation*}
Splitting the integration domain and using the identity $\cos\theta-1 = -\sin^2\theta/(1+\cos\theta)$ we obtain
\begin{equation*}
\begin{split}
 I_3\left[u_3^{\downarrow};S_3^R\right](\ner)\sim&\ \sqrt{\frac{k_3R}{8\pi}}\e^{ik_3R-i\frac\pi4}\left\{-\int_{-\frac\pi2}^{0}\frac{\sin^2\theta_3}{1+\cos\theta_3}\e^{-ik_3(d_{2}\sin\theta_3-\ner\cdot\hat\ner')+iRk_3\cos\theta_3}\de\theta_3\right.+\\
&\left.\int_{-\pi}^{-\frac\pi2}\left\{\cos\theta_3-1\right\}\e^{-ik_3(d_{2}\sin\theta_3+\ner\cdot\hat\ner')+iRk_3\cos\theta_3}\de\theta_3\right\}.
\end{split}
\end{equation*}
 Integration by parts yields that  the first  integral above amounts to a quantity of order $ O\lf(R^{-1/2}\rg)$. The stationary point at $\theta=-\pi$ in the second integral, on the other hand, leads to
\begin{equation*}
 I_3\left[u_3^{\downarrow};S_3^R\right](\ner)=-\frac{\e^{ik_3x}}{2}+ O\lf(R^{-1/2}\rg).
\end{equation*}
Similarly, in the case $k_3=-k_1\cos\alpha$ for
$\alpha\in(-\pi,\pi/2]$ it can be shown that $
I_3\left(u_3^{\downarrow};S_3^R\right)=-\e^{-ik_3x}/2+
O\lf(R^{-1/2}\rg)$. In the case~(c), finally, $u_3^{\downarrow}(\ner')
=a\e^{ik_3R\cos(\theta_3-\alpha')}$, $\ner'\in S_3^R$, where
$a=\e^{-ik_3d_{2}\sin\alpha'}$ and where the angle
$\alpha'\in(-\pi,0)$ is determined by the Snell's law $k_3\cos\alpha'
=k_1\cos\alpha $. Thus, once again, integration by parts yields
\begin{align*}
I_3\lf[u_3^{\downarrow};S^R_3\rg](\ner) \sim&\ a\sqrt{\frac{k_3R}{8\pi}}\e^{ik_3R-i\frac\pi4}\!\int\limits_{-\pi}^{0}\left\{\cos(\theta_3-\alpha')-1\right\}\e^{-ik_3(d_{2}\sin\theta_3 +  \ner\cdot\hat\ner')} \e^{ik_3R\cos(\theta_3-\alpha')}\!\de\theta_3\\
=&\  O\left(R^{-1/2}\right).
\end{align*}

Therefore, taking the limit as $R\to\infty$
in~\eqref{eq:green_truncated} we obtain
\begin{subequations}\begin{eqnarray}
 I_1\left[v_1;\Gamma_1\right](\ner) +\delta u^\inc(\ner)&=&\left\{\begin{array}{cl} v_1(\ner), &\ner\in\Omega_1,\smallskip\\
 0,& \ner\in\R^2\setminus\overline\Omega_1,
 \end{array}\right.\label{eq:green_1_3L}\\
 I_2\left[v_2;\Gamma_{2}\right](\ner)- I_2\left[v_2;\Gamma_{1}\right](\ner)  &=&\left\{\begin{array}{cl} v_2(\ner), &\ner\in\Omega_2,\smallskip\\
 0,& \ner\in\R^2\setminus\overline\Omega_2,
 \end{array}\right.\label{eq:green_j_3L}\\
 -I_3\left[v_3;\Gamma_{2}\right](\ner)+\delta u^{\parallel}_3(\ner) &=&\left\{\begin{array}{cl} v_3(\ner), &\ner\in\Omega_3,\smallskip\\
 0,& \ner\in\R^2\setminus\overline\Omega_3,
 \end{array}\right.\label{eq:green_N_3L}
\end{eqnarray}\label{eq:green_representation}\end{subequations}
where $u^\parallel_3$ in~\eqref{eq:green_N_3L} is given by 
\begin{equation}
u^{\parallel}_N(\ner) = \left\{\begin{array}{ccc}
\displaystyle\frac{q_N\e^{ik_{1}x\cos\alpha}}{2}&\mbox{if}&k_N=k_1|\cos\alpha|\medskip,\\
0&\mbox{if}&k_N\neq k_1|\cos\alpha|,
\end{array}\right.\label{eq:parallel_ML}
\end{equation}
with $N=3$, and where $\delta=0$ if $v_j= v_j^d$, and $\delta=1$ if
$v_j=\tilde u_j^p+v_j^d$ or $v_j=\tilde u_j^p$.
\begin{remark}
  The total field representation~\eqref{eq:green_representation} can
  easily be extended to problems of scattering by defects in the
  presence of layer media composed by $N>3$ layers; the result is
\begin{subequations}\begin{eqnarray}
 I_1\left[v_1;\Gamma_1\right](\ner) +\delta u^\inc(\ner)&=&\left\{\begin{array}{cl} v_1(\ner), &\ner\in\Omega_1,\smallskip\\
 0,& \ner\in\R^2\setminus\overline\Omega_1,
 \end{array}\right.\label{eq:green_1_NL}\\
 I_{j}\left[v_j;\Gamma_{j}\right](\ner)- I_j\left[v_j;\Gamma_{j-1}\right](\ner)  &=&\left\{\begin{array}{cl} v_j(\ner), &\ner\in\Omega_j,\smallskip\\
 0,& \ner\in\R^2\setminus\overline\Omega_j,
 \end{array}\right. \  j=2,\ldots,N-1,\label{eq:green_j_NL}\\
 -I_N\left[v_N;\Gamma_{N-1}\right](\ner)+ \delta u^{\parallel}_N(\ner) &=&\left\{\begin{array}{cl} v_N(\ner), &\ner\in\Omega_N,\smallskip\\
 0,& \ner\in\R^2\setminus\overline\Omega_N.
 \end{array}\right.\label{eq:green_N_NL}
\end{eqnarray}\label{eq:green_representation_NL}\end{subequations}

\end{remark}

\bibliographystyle{abbrv}
\bibliography{References}

\end{document}